\documentclass[11pt]{article}

\usepackage{amsmath,amsfonts,amsbsy,amstext,amsopn,amssymb}

\usepackage{longtable}
\usepackage{theorem}
\usepackage{algorithmic}
\usepackage{algorithm}

\makeatletter \@addtoreset{equation}{section}

\makeatother

{\theorembodyfont{\slshape}
\newtheorem{theorem}{Theorem}
\newtheorem{proposition}{Proposition}
\newtheorem{lemma}{Lemma}

}
{\theorembodyfont{\rmfamily}

\newtheorem{remark}{Remark}
\newtheorem{example}{Example}

}

\newcommand{\vect}{{\sf{vec}}}
\newcommand{\tr}{{\sf{trace}}}

\newcommand{\rt}{{\top }}

\def\e{\varepsilon}
\def\Oh{{\mathcal O}}

\def\Prob{\mathop{\rm Prob}}
\def\Dom{\mbox{\sf Dom}}

\def\bfe{{\mathbf e}}

\def\hankel{{\sf Hankel}}

\def\VdM{{\sf VdM}}
\def\Cauchy{{\sf Cauchy}}

\def\SS{\mathcal S}
\def\LL{\mathcal L}
\def\str{{\LL}}

\def\Diag{\mbox{\sf Diag}}

\def\bfe{{\mathbf e}}

\def\reg{{\sf Reg}}

\def\sfv{{\sf v}}

\usepackage{multirow}

\newcommand {\eq} [1] {\begin{equation}\label{#1}}
\newcommand {\en} {\end{equation}}
\newcommand {\proof} {\par{\it Proof}. \ignorespaces}
\newcommand {\eproof}
      {\space
        {\ \vbox{\hrule\hbox{\vrule height1.3ex\hskip0.8ex\vrule}\hrule}}
        \par}
\newcommand {\R}        {{\mathbb R}}

\newcommand {\Rm}       {\R^m}

\newcommand {\Rpn}      {\R^{p \times n}}

\newcommand {\Rpq}      {\R^{p \times q}}

\newcommand {\Rmn}      {\R^{m \times n}}

\newcommand {\mat}      [1] {\left[\begin{array}{#1}}
\newcommand {\rix}          {\end{array}\right]}

\newcommand {\rank}     {\mathop{\rm rank}\nolimits}

 \catcode`@=11
 \font\tenex=cmex10 
 \newdimen\p@renwd
 \setbox0=\hbox{\tenex B} \p@renwd=\wd0 
 \def\bmat#1{\begingroup \m@th
   \setbox\z@\vbox{\def\cr{\crcr\noalign{\kern2\p@\global\let\cr\endline}}%
     \ialign{$##$\hfil\kern2\p@\kern\p@renwd&\thinspace\hfil$##$\hfil
       &&\quad\hfil$##$\hfil\crcr
       \omit\strut\hfil\crcr\noalign{\kern-\baselineskip}%
       #1\crcr\omit\strut\cr}}%
   \setbox\tw@\vbox{\unvcopy\z@\global\setbox\@ne\lastbox}%
   \setbox\tw@\hbox{\unhbox\@ne\unskip\global\setbox\@ne\lastbox}%
   \setbox\tw@\hbox{$\kern\wd\@ne\kern-\p@renwd\left[\kern-\wd\@ne
     \global\setbox\@ne\vbox{\box\@ne\kern2\p@}%
     \vcenter{\kern-\ht\@ne\unvbox\z@\kern-\baselineskip}\,\right]$}%
   \null\;\vbox{\kern\ht\@ne\box\tw@}\endgroup}
 \catcode`@=12
\begin{document}

\title{\bf  Structured Condition Numbers
of Structured Tikhonov Regularization Problem and their Estimations}
\author{Huai-An Diao
\thanks{School of Mathematics and Statistics \& Key Laboratory
for Applied Statistics of MOE, Northeast
Normal University, No. 5268 Renmin Street, Chang Chun 130024,
P. R. of China. ({\tt hadiao@nenu.edu.cn and hadiao78@yahoo.com})
This author is supported by the
National Natural Science Foundation of China. Part of this work was
finished when the author visited Shanghai
Key Laboratory of Contemporary Applied Mathematics and
McMaster University in  2014.}
\and
Yimin Wei
\thanks{School of Mathematics
\& Shanghai Key Laboratory of Contemporary Applied Mathematics,
Fudan University, Shanghai 200433, P. R. of China. ({\tt
ymwei@fudan.edu.cn and yimin.wei@gmail.com}).  This author is
supported by the National Natural Science Foundation of China under
grant 11271084. }
\and
Sanzheng Qiao\thanks{
Department of Computing and Software,
McMaster University, Hamilton, Ontario L8S4K1, Canada.
({\tt qiao@mcmaster.ca}) This author is partially supported by
Natural Science and Engineering Council (NSERC) of Canada
and Shanghai Key Laboratory of Contemporary Applied Mathematics.}
}
\date{}
\maketitle

\begin{quote}
{\bf Abstract.} Both structured componentwise and structured normwise
perturbation analysis of the Tikhonov regularization  are presented.
The structured matrices under consideration include:
Toeplitz, Hankel, Vandermonde, and
Cauchy matrices. Structured normwise, mixed and componentwise
condition numbers for the Tikhonov regularization are introduced and
their explicit expressions are derived. For the general linear structure,
we prove the structured condition numbers are smaller than
their corresponding unstructured counterparts based on the
derived expressions. By means of the power method and
small sample condition estimation,
the fast condition estimation algorithms are proposed.
Our estimation methods can
be integrated into Tikhonov regularization algorithms that use the
generalized singular value decomposition
(GSVD). The structured condition numbers and
perturbation bounds are tested on some numerical examples and
compared with their unstructured counterparts. Our numerical examples
demonstrate that the structured mixed condition
numbers give sharper perturbation bounds than existing ones,
and the proposed condition estimation algorithms are reliable.
\end{quote}

{\small {\bf Keywords:}
Tikhonov regularization, structured matrix, condition number,
componentwise, structured perturbation, small sample condition estimation.\\

{\bf AMS Subject Classification:} 15A09, 15A12, 65F35.}


\section{Introduction}
\setcounter{equation}{0}

For discrete ill-posed problems, the Tikhonov regularization
(cf. \cite{Tikhonov}) reads
\begin{equation}\label{TikReg1}
\min\limits_{x}\left\{\|Ax-b\|_2^2+\lambda^2\|Lx\|_2^2\right\},
\qquad A\in\R^{m\times n} \quad {\rm and } \quad  L\in\R^{p\times n}
\end{equation}
where $\lambda$  is the regularization parameter,
which controls  the weight between $\|Lx\|_2$ and the residual $\|Ax-b\|_2$.
The matrix $L$ is typically the identity
matrix $I_n$ or a discrete approximation to some derivation
operator. Tikhonov regularization is also known as {\em ridge
regression} in statistics~\cite{BjorckBook}.

For the regularization problem (\ref{TikReg1}),
to ensure the uniqueness of the solution for
any $\lambda > 0$, we
always assume that $\rank(L)=p\leq n\leq m$ and $\rank\left(\mat{c} A \\
L \rix\right) =n$ (cf. \cite[\S5]{BjorckBook}).
The regularization problem (\ref{TikReg1}) can be
rewritten in the matrix form
\begin{equation}\label{TikReg2}
\min\limits_{x}\left\| \mat{c} A\\ \lambda L \rix x-\mat{c} b\\
{\bf 0}\rix \right\|_2,
\end{equation}
where $\bf0$ is the zero vector. Since the normal equations
corresponding to (\ref{TikReg2}) are
\begin{equation}\label{NormEqn}
\left(A^\rt A+\lambda^2L^\rt L \right)x=A^\rt b,
\end{equation}
we can obtain the following explicit expression for the {\em
Tikhonov regularized solution}:
$$
x_{\lambda}=\left(A^\rt A+\lambda^2 L^\rt L \right)^{-1}A^\rt b.
$$

Alternatively, the problem (\ref{TikReg1}) can also be solved by the
generalized singular value decomposition
(GSVD)~\cite{hansengsvd,Hansen,VanLoan}. For  rectangular matrices
$A \in \R^{m\times n}$ and $L \in \R^{p \times n}$ with
$\rank(L)=p$ and $\rank\left(\mat{c} A \\
L \rix\right) =n$, the GSVD of $(A,L)$ is given by the pair of factorizations
\begin{equation}\label{GSVD}
A=U\begin{bmatrix}\Sigma& {\bf 0}\cr {\bf 0}& I_{n-p}
\end{bmatrix}RQ^\rt \quad \hbox{and} \quad
L=V\begin{bmatrix}S &{\bf
0}\end{bmatrix}RQ^\rt,
\end{equation}
where $U \in \R^{m\times n}$ has orthonormal columns,
$V \in \R^{p \times p}$, $Q\in \R^{n \times n} $ are orthogonal, $R$ is $n$-by-$n$,
upper triangular and nonsingular, and $\Sigma$ and $S$ are
$p\times p$ diagonal matrices:
$\Sigma=\Diag(\sigma_1,\sigma_2,\ldots,\sigma_{p})$ and
$S=\Diag(\mu_1,\mu_2,\ldots,\mu_p)$ with
$$
0 \leq \sigma_1 \leq \sigma_2 \leq \ldots \leq \sigma_p < 1
\quad \hbox{and} \quad
 1 \geq \mu_1 \geq \mu_2 \geq \ldots \geq \mu_p>0,
$$
satisfying $\Sigma^2+S^2=I_{p}$. Then the
{\em generalized singular values $\gamma_i$} of $(A,L)$ are defined
by the ratios $\gamma_i=\sigma_i/\mu_i$ ($i=1,2,\ldots,p$).  Once
the GSVD is computed, the Tikhonov regularized solution can be obtained
by \cite[Chapter 4]{Hansen}
$$
x_\lambda=QR^{-1}\begin{bmatrix}F&{\bf 0}\cr {\bf 0}
&I_{n-p}\end{bmatrix}
\begin{bmatrix}\Sigma^{\dagger}&{\bf 0}\cr {\bf 0}&I_{n-p}\end{bmatrix}U^\rt
b, \quad F=\Diag(f_1,f_2,\ldots,f_p),
$$
where $f_i= \gamma_i^2 / ( \gamma_i^2+\lambda^2 )$
for $i=1,2,\ldots, p$, are called the {\em filter factors} for
the Tikhonov regularization \cite{Hansen,Hansen2} and
$\Sigma^\dagger$ is the
Moore-Penrose inverse of $\Sigma$ \cite{BjorckBook}.

In sensitivity analysis, condition numbers are of great importance
because they measure the {\em worst-case} effect of small changes
in the data on the solution. For the perturbation analysis of
the linear least squares (LS) problem, the
reader is referred to~\cite{Arioli,Gratton,Baboulin09,CD07,CDW07}.
Arioli et al.~\cite{Arioli} introduced a {\em partial condition number}
of the LS
problem, which can be viewed as a condition number of a {\em linear
functional} of the LS problem. Baboulin et al.~\cite{Gratton} have shown
that the partial condition numbers of the LS problem represent some
quantities in statistics.  For the perturbation analysis for the Tikhonov
regularization, we refer to \cite{Marten,Inverse} and references
therein. Malyshev~\cite{Malyshev} adopted a unified
theory to study the normwise condition numbers for the Tikhonov
regularization. Chu et al.~\cite{CLTW}
investigated the componentwise perturbation analysis of the
Tikhonov regularization problems and derived condition
number expressions involving the Kronecker products, which can be of huge
dimension even for small problems, preventing us from
estimating the condition numbers while solving the Tikhonov
regularization problem.  In this paper, we consider the structured
condition numbers for a linear functional of the Tikhonov
regularization. Fast condition number estimation, which is important
in practice, is discussed.

Structured matrix computation is a hot research topic; see \cite{Jin,Ng} and the references therein. The structured Tikhonov regularization problem was recently
studied in~\cite{Benzi,Chan,hansentoep}. Eld\'{e}n gave a stable efficient algorithm for the Tikhonov
regularization with triangular Toeplitz structure. Park and Eld\'{e}n \cite{Park00}  devised fast algorithms for solving LS with Toeplitz structure, based on the generalization of the classical Schur algorithm,  and discussed their stability properties. Also, Park and Eld\'{e}n studied the stability analysis and fast algorithms
for triangularization of rectangular Toeplitz matrices \cite{Park97}. Hence, it is natural
to investigate {\em structured perturbations} on the structured coefficient
matrix, which lead to the structured condition numbers for
the structured Tikhonov regularization problem. Structured condition
numbers for several categories  of structured matrices have been
presented in
\cite{Desmond,BG,BF,CD07,Gohberg,HighamHigham,Rump03a,Rump03b,Hua,WWQ}.
In this paper we derive explicit formulas for the condition numbers of
the Tikhonov regularization problem, when perturbations of $(A,b)$
are measured by normwise or componentwise or a mixture of normwise and
componentwise. To make our discussion general, we consider the condition
number of $Mx$, i.e., a linear function of the Tikhonov regularized
solution, where $M \in \R^{l\times n}$ and $x\in \R^{n}$, $l\leq
n$. The common situations are the special cases, when $M$
is the identity matrix (condition number of the Tikhonov regularized
solution) or a canonical vector (condition number of one component of
the solution). We obtain the expressions of the structured condition
numbers in the absence of the Kronecker product, so that they can be
estimated by the power method due to Hager~\cite{Hager} and
Higham~\cite{Higham_ACM,Higham_SISC},  see~\cite[Chapter
15]{HighamBook} for the detail, while solving the Tikhonov
regularization problem.

Moreover, in this paper, we adopt the statistical condition
estimation (SCE) method \cite{KenneyLaub_SISC94} for
numerically estimating the condition of Tikhonov regularization problem.
The SCE can be used to estimate the componentwise local sensitivity of
any differentiable function at a given input data, which is flexible and
accommodates a wide range of perturbation types such as
structured perturbations.
Thus SCE often provides less conservative estimates than the methods that
do not exploit structures. The SCE method has been shown to
be both reliable and efficient for many problems including linear systems
\cite{KenneyLinear}, structured linear systems \cite{laublinear}, linear
least squares problems \cite{Kenney}, eigenvalue problems
\cite{Gudmundsson,laubeig}, matrix functions \cite{KenneyLaub_SISC94},
the roots of polynomials \cite{laubroot}, etc.

We follow the convention of representing a point
$x\in\R^n$ as a column vector. If $x\in\R^n$ and $y\in\R^m$, then
 $[x;y]$ is an $m+n$ column vector by stacking $x$ on top of $y$.
If $A\in\R^{m\times n}$ and $B \in\R^{m\times
q}$, then $[A,B]$ denotes the matrix obtained by putting $A$ and $B$
side by side. The symbol  `$.^\rt$' denotes matrix transpose, $\|\cdot\|_2$
is the spectral norm, $\|\cdot\|_F$
is the Frobenius norm and $\|\cdot\|_\infty$ is the infinity norm. The matrix
$\Diag(d)\in\R^{q\times q}$ denotes a diagonal matrix with
the vector $d$'s entries being its corresponding diagonal components.
For any points $a,b
\in \R^n$, the vector $c=\frac{a}{b}$ is obtained by
componentwise division. In particular, $b_i=0$ assumes $a_i=0$,
and in this case $c_i=0$.
For a matrix $A\in\R^{m\times n}$, we define
$\vect(A)\in\R^{mn}$ by $\vect(A)=[a_1^\rt,
a_2^\rt,\ldots,a_n^\rt]^\rt$, where $A=[a_1,a_2,\ldots, a_n]$ with
$a_i \in \R^m$, $i=1,2,\ldots,n$. The ${\sf unvec}$ operation
is defined as $A={{\sf unvec}}(v)$ which sets the entries of $A$ to
$a_{ij}=v_{i+(j-1)n}$ for
$v=[v_1,v_2,\ldots,v_{mn}] \in \R^{1 \times mn}$.
We define a permutation matrix $\Pi$ of order $mn$ so that
$\Pi(\vect(A))=\vect\left(A^\rt\right)$.
Let `$\otimes$' denote the {\em Kronecker product}~\cite{Graham}, i.e.,
$A\otimes B=[a_{ij}B]\in\R^{mp\times nq}$ for $A=(a_{ij})\in
\R^{m\times n}$ and $ B \in \R^{p \times q}.$ The notation $|A|\leq
|B|$  means that $|a_{ij}|\leq |b_{ij}|$. For the Kronecker product,
we recall the following properties, which can
be found in~\cite{Graham},
\begin{equation}\label{eq:kr}
(A\otimes B)^\rt=A^\rt
\otimes B^\rt, \quad |A \otimes B| = |A| \otimes |B|,
\quad \vect (AXB) = \left(B^\rt \otimes A\right) \vect
(X),
\end{equation}
where $|A|=[|a_{ij}|]$ and $a_{ij}$ is the $(i,j)$-th entry of $A$.

This paper is organized as follows. We provide some preliminaries in
Section \ref{sec:preliminaries}, investigate matrices with linear structures
in Section \ref{linear structure: section title} and move to matrices with
nonlinear structures in Section \ref{sec:nonlinear}.
The SCE-based condition estimation algorithms are proposed in
Section \ref{sec:SCE for Tik}. In Section \ref{sec:numerical exam},
we demonstrate test results showing the sharpness of our
structured condition numbers and effectiveness of the condition
estimation algorithms. Finally, conclusions are drawn in the last section.

\section{Preliminaries}\label{sec:preliminaries}

In this section, we first recall the general (unstructured)
condition number definitions \cite{Gohberg}.
Then we consider the structured Tikhonov regularization problems,
introduce structured perturbations, and define their structured
condition numbers. Finally, we briefly describe the basic ideas of SCE.

\subsection{Structured condition numbers for the Tikhonov regularization}

For $x,\, a \in \R^p$ and $\e>0$ we denote\
$S(a,\e)=\{x\in \R^p \mid |x-a|\leq \e |a|\}$ and
$T(a,\e)=\{x \in \R^p ~| ~\|x-a\|_2 \leq \e\}$. For a function
$F:\R^p \rightarrow \R^q$, we denote $\Dom(F)$ as its domain.
The following lemma defines general (unstructured) condition numbers.

\begin{lemma}\label{def:CN}{\rm \bf (\cite{Gohberg})}
Let $F:\R^p \rightarrow \R^q$ be a continuous mapping defined on an
open set $\Dom(F) \subset \R^p$. Let $a \in \Dom(F)$ such that $a
\neq \bf 0$ and $F(a)\neq \bf 0$.
\begin{description}
\item[(i)] The {\em mixed condition number} of $F$ at $a$ is
defined by
$$
m(F,a)=\lim_{\e \rightarrow 0} \sup_{x \in S(a,\e) \atop x \neq a}
\frac{\|F(x)-F(a)\|_\infty}{\|F(a)\|_\infty}\frac{1}{d(x,a)}=
\frac{\||{\bf D}{F}(a)|~|a|\|_\infty}{\|F(a)\|_\infty},
$$
where ${\bf D}{F}(a)$ is the Fr\'{e}chet derivative of $F$ at $a$
and $|a|=(|a_i|)$ with $a=[a_1,a_2,\ldots,a_p]^\rt$.
\item[(ii)]
Suppose $F(a)=(f_1(a),f_2(a),\ldots,f_q(a))$ is such that
$f_j(a)\neq 0$ for $j=1,2,\ldots,q$. Then the {\em componentwise
condition number} of $F$ at $a$ is
$$
c(F,a)=\lim_{\e \rightarrow 0} \sup_{x \in S(a,\e) \atop x \neq a}
\frac{d(F(x),F(a))}{d(x,a)}=
\left\|\frac{|{\bf D}F(a)||a|}{|F(a)|}\right\|_\infty.
$$
\item[(iii)]
The {\em normwise condition number} of $F$ at $a$ is defined by
$$
\kappa(F,a)=\lim_{\e \rightarrow 0} \sup_{x \in T(a,\e) \atop x \neq a}
\frac{\|F(x)-F(a)\|_2}{\|x-a\|_2}\frac{\|a\|_2}{\|F(a)\|_2}=
\frac{\|{\bf D}{F}(a)\|_2\|a\|_2}{\|F(a)\|_2}.
$$
\end{description}
\end{lemma}

In the following we assume that $\Delta A$ and $\Delta b$ are
perturbations to  $A$ and $b$ respectively, which satisfy
$\rank \left(\mat{c} A+\Delta A \\ L \rix\right)=n$.
The perturbed counterpart of the problem (\ref{TikReg1}) and its
normal equations (\ref{NormEqn}) are, respectively,
\begin{equation}\label{TikReg1-1}
\min\limits_{x+\Delta x}\left\{\|(A+\Delta A)(x+ \Delta x)-
(b+\Delta b)\|_2^2+\lambda^2\|L(x+\Delta x)\|_2^2\right\},
\end{equation}
and
\[
\Big[(A+\Delta A)^\rt (A+\Delta A)+\lambda^2L^\rt
L\Big](x_{\lambda}+\Delta x )=(A+\Delta A)^\rt (b+\Delta b).
\]
Then the perturbed Tikhonov regularized solution is given by
\begin{equation}\label{deltax}
x_\lambda+\Delta x=\left[(A+\Delta A)^\rt (A+\Delta
A)+\lambda^2L^\rt L \right] ^{-1}(A+\Delta A)^\rt (b+\Delta b).
\end{equation}

Denoting
\[
P(A,\lambda)=\left(A^\rt A+\lambda^2 L^\rt L \right)^{-1},
\]
Chu~et al. \cite{CLTW} define the {\em non-structured}
mixed, componentwise, and normwise condition
numbers for the Tikhonov regularization and obtain respectively
\begin{align}
{\rm m_{Reg}} &=\lim_{\epsilon \rightarrow 0}\sup_{|\Delta
A|\leq\epsilon|A|\atop{|\Delta b|\leq\epsilon|b|}}
\frac{\|\Delta x\|_\infty}{\epsilon\|x_{\lambda}\|_\infty}=
\frac{\left\|~|H(A,b)|\vect(|A|)+\left|P(A,\lambda)A^\rt
\right|\,|b|~\right\|_\infty}
        {\|x_{\lambda}\|_\infty},\label{mReg}\\
{\rm c_{Reg}}&=\lim_{\epsilon \rightarrow 0}\sup_{|\Delta A|
\leq\epsilon|A|\atop{|\Delta
b|\leq\epsilon|b|}}\frac{1}{\epsilon}\left\|\frac{\Delta x}
{x_{\lambda}}\right\|_\infty=\left\|\frac{|H(A,b)|
\vect(|A|)+\left |P(A,\lambda)A^\rt \right||b|}
{x_\lambda}\right\|_\infty, \label{cReg}\\
{\rm  cond^F_{Reg}} &=\lim_{\epsilon\rightarrow0}
\sup_{\|\left[\Delta A\,,\, \Delta b\right]\|_F\leq
\epsilon\|\left[A \,,\, b\right]\|_F}\frac{\|\Delta x\|_2}
{\epsilon\|x_{\lambda}\|_2}=\frac{\left \| \left[ H(A,b),\,
P(A,\lambda)A^\rt \right]\right \|_2\left \|[ A ,\, b ]\right \|_F}
{\|x_{\lambda}\|_2}, \label{condF}
\end{align}
where $H(A,b)=-x_\lambda^\rt \otimes \left[P(A,
\lambda)A^\rt\right]+\left[ P(A, \lambda)\otimes
r_\lambda^\rt\right]$ and $\,r_\lambda=b-Ax_\lambda$.

If we define a mapping
\begin{equation}\label{psi}
\psi:\  [A,\, b] \in \R^{m\times n}\times\R^m \mapsto
\left(A^\rt A+\lambda^2 L^\rt L
\right)^{-1}A^\rt b \in \R^n
\end{equation}
then it is easy to see that the definitions in Lemma~\ref{def:CN} are
equivalent to \eqref{mReg}-\eqref{condF}, that is,
$${\rm m_{Reg}}:=m(\psi,[A,\,b]), \quad
{\rm c_{Reg}}:=c(\psi,[A,\,b]), \quad
{\rm cond^F_{Reg}}:=\kappa(\psi,[A,\,b]).$$

When the coefficient matrix $A$ in~\eqref{TikReg1} has some
structures, such as Toeplitz, it is reasonable to assume that the
perturbation $\Delta A$ in~\eqref{TikReg1-1} has the same structure
of $A$. Then $\Delta A$ is called {\it structured perturbation}
\cite{Rump03a,Rump03b} on $A$. Usually a structured matrix
$A \in \mathbb{R}^{m \times n}$ can
be represented by fewer than $mn$ parameters. For example,
an $m \times n$ Toeplitz matrix can be represented by its
first column and last row, $m+n-1$ parameters.
Here we use a mapping to characterize this
relationship. Let $\SS$ be the set of structured matrices
under consideration and $a$ the vector representing
a structured matrix $A$, then we define a mapping
$$
g: \  a \in \R^k \mapsto A \in \SS .
$$

In order to apply Lemma~\ref{def:CN} to define the structured
condition numbers for the Tikhonov regularization,
we construct a mapping
\begin{equation}\label{phi}
\phi: \  [a;b] \in \R^{k+m} \mapsto
M \left(A^\rt A+\lambda^2 L^\rt
L\right)^{-1}A^\rt b \in \R^l,
\end{equation}
where $M\in\R^{l\times n}$, $l \leq n$, is general. In particular,
when $M=e_i^{\rm T}$, the $i$-th column of the identity matrix, then
we are interested in some particular component of $x_\lambda$.

Let $\Delta a$ be the perturbation on $a$, then the structured
perturbation matrix $\Delta A$ on $A$ in \eqref{TikReg1-1} is
$g(a+\Delta a)-g(a)$. Now we are ready to define the {\em structured mixed,
componentwise and normwise} condition numbers for a linear
functional of the structured Tikhonov regularization,

\begin{align*}
m^{\reg}_\SS(A,b)&:=m(\phi,[a;b])=\lim_{\epsilon \rightarrow
0}\sup_{|\Delta a|\leq\epsilon|a|\atop{|\Delta
b|\leq\epsilon|b|}}\frac{\left\|M\Delta
x\right\|_\infty}{\epsilon \left\|Mx_{\lambda}\right\|_\infty},\cr
c^{\reg}_\SS(A,b)&:=c(\phi,[a;b])=\lim_{\epsilon \rightarrow
0}\sup_{|\Delta a|\leq\epsilon|a|\atop{|\Delta
b|\leq\epsilon|b|}}\frac{1}{\epsilon}\left\|\frac{M\Delta
x}{Mx_{\lambda}}\right\|_\infty,\cr
\kappa^{\reg}_\SS(A,b)&:=\kappa(\phi,[a;b])
=\lim_{\epsilon\rightarrow0}\sup_{\|\left[\Delta a; \,\Delta
b\right]\|_2\leq\epsilon\|\left[a ;\, b\right]\|_2}\frac{\|M\Delta
x\|_2}{\epsilon\|Mx_{\lambda}\|_2},
\end{align*}
where $\Delta x$ is defined in \eqref{deltax}.

\begin{remark}\label{remark:1}
Note that here $g$ is a general mapping, in that it can represent
any structure. When  the structure in $A$ is linear, such as
symmetric, or Teoplitz, or Hankel, we can choose $g$ a linear mapping,
which will be discussed in Section \ref{linear structure: section title}.
When $A$ has a nonlinear structure
such as Vandermonde or Cauchy, we can choose a nonlinear mapping $g$
to define the structured condition numbers. Especially we can define
the unstructured linear functional condition number for $x_\lambda$
when we restrict $\SS$ to be $\R^{m\times n}$,
which are generalizations of \eqref{mReg}, \eqref{cReg}
and \eqref{condF}, as follows
\begin{align*}
m^{\reg}(A,b)&=\lim_{\epsilon \rightarrow
0}\sup_{|\Delta A|\leq\epsilon|A|\atop{|\Delta
b|\leq\epsilon|b|}}\frac{1}{\epsilon}\frac{\left\|M\Delta
x\right\|_\infty}{\left\|Mx_{\lambda}\right\|_\infty},\quad
c^{\reg}(A,b)=\lim_{\epsilon \rightarrow
0}\sup_{|\Delta A|\leq\epsilon|A|\atop{|\Delta
b|\leq\epsilon|b|}}\frac{1}{\epsilon}\left\|\frac{M\Delta
x}{Mx_{\lambda}}\right\|_\infty,\cr
\kappa^{\reg}(A,b)
&=\lim_{\epsilon\rightarrow0}\sup_{\|\left[\Delta A ,\Delta
b\right]\|_2\leq\epsilon\|\left[A , b\right]\|_2}\frac{\|M\Delta
x\|_2}{\epsilon\|Mx_{\lambda}\|_2}.
\end{align*}
When $M=I_n$, the above definitions reduce to \eqref{mReg},
\eqref{cReg}  and \eqref{condF}.
\end{remark}

Finally, we give the well-known Banach lemma,
which will be useful in  Section \ref{linear structure: section title}.
\begin{lemma}\label{lemma:Banach}
Let $E\in \R^{n\times n}$ and $\|\cdot\|$ be any norm on $\R^{n\times n}$,
if $\|E\|<1$, then
$I_n+E$ is nonsingular and its inverse can be expressed by
$$
\left(I_n+E\right)^{-1}=I_n-E+\Oh(\|E\|^2).
$$
\end{lemma}

\subsection{Statistical condition estimation}\label{sec:SCE}

In SCE, a small random perturbation is introduced to the input, and
the change in the output, by an appropriate scaling, is measured as
a condition estimate. Explicit bounds on the probability of
the accuracy of the estimate exist \cite{KenneyLaub_SISC94}.
The idea of SCE can be illustrated by a general real-valued function:
$f:\R^{p}\rightarrow \R$, and we are interested in the sensitivity at
some input vector $x$. By the Taylor theorem we have
$$
f(x+\delta d)-f(x)=\delta ({\bf D} f(x))^\rt d+\Oh(\delta ^2),
$$
where $\delta$ is a small scalar, $\|d\|_2=1$ and ${\bf D} f(x)$ is
the  Fr\'{e}chet derivative of $f$ at
$x$. Note that the quantity $({\bf D} f(x))^\rt d$
(denoted by ${\bf D} f(x;d) $) is
just the directional derivative of $f$ with respect to $x$ at
the direction $d$.
It is easy to see that up to the first order in $\delta $,
$$
\left|f(x+\delta d)-f(x)\right|\approx \delta {\bf D} f(x;d),
$$
then the local sensitivity can be measured by
$\|{\bf D} f(x)\|_2$. The condition numbers of $f$ at $x$ are mainly
determined by the norm of the gradient ${\bf D} f(x)$
(\cite{KenneyLaub_SISC94}). According to~\cite{KenneyLaub_SISC94},
if we select $d$ uniformly and randomly from the unit $p$-sphere
$S_{p-1}$ (denoted $d \in {\cal U}(S_{p-1})$), then the expectation
${\bf E}(|{\bf D} f(x;d)|/\omega_p)$  is $\|{\bf D} f(x)\|_2$,
where $\omega_p$ is the Wallis factor.
In practice, the Wallis factor can be
approximated accurately \cite{KenneyLaub_SISC94} by
\[
\omega_p\approx
\sqrt{\frac{2}{\pi(p-\frac{1}{2})}}.
\]
Therefore,  we can use
$$
\nu=\frac{|{\bf D} f(x;d)|}{\omega_p}
$$
as a condition estimator, which can estimate $\|{\bf D} f(x)\|_2$
with high probability for the function $f$ at $x$ (see
\cite{KenneyLaub_SISC94} for details), for example,
$$
\Prob\left(\frac{\|{\bf D} f(x)\|_2}{\gamma}\leq \nu \leq \gamma
\|{\bf D}f(x)\|_2\right) \geq 1-\frac{2}{\pi
\gamma}+\Oh\left(\frac{1}{\gamma^2}\right) ,
$$
for $\gamma > 1$.
We can use multiple samples of $d$, denoted $d_j$, to increase the accuracy
\cite{KenneyLaub_SISC94}. The $t$-sample condition estimation is given by
$$
\nu(k)=\frac{\omega_t}{\omega_p}\sqrt{|{\bf D} f(x;d_1)|^2
+|{\bf D} f(x;d_2)|^2 + \cdots + |{\bf D} f(x;d_t)|^2},
$$
where $[d_1,d_2,\ldots,d_t]$ is orthonormalized after
$d_1,d_2,\ldots,d_t$ are selected uniformly and
randomly from ${\cal U}(S_{p-1})$. The accuracy of $\nu(2)$ is given by
\begin{align*}
\Prob\left(\frac{\|\nabla f(x)\|_2}{\gamma}\leq \nu(2) \leq \gamma
\|\nabla f(x)\|_2\right) &\approx 1-\frac{\pi}{4 \gamma^2},\quad \gamma>1.
\end{align*}
Usually, a few samples are sufficient for good accuracy. These results can be conveniently
generalized to vector- or matrix-valued  functions by viewing $f$ as a map
from $\R^p$ to $\R^q$. The operations {\sf vec} and {\sf unvec} can be used
to convert between matrices and
vectors, where each of the $q$ entries of $f$ is a scalar-valued function. Evaluating the
matrix function at a slightly perturbed argument yields a local condition estimate for
one component of the computed solution.

\section{Linear Structures}\label{linear structure: section title}

In this section, we consider the classes $\LL$ of structured matrices
that is a linear subspace of $\R^{m\times n}$. The examples of such class
include Toeplitz  and Hankel matrices. We first present a structured
perturbation analysis and structured condition numbers. Then we
propose efficient condition number estimators using the power method.

\subsection{Condition numbers}

Suppose that $\dim(\LL)=k$, and $S_{1},S_{2},\ldots,S_{k}$
form a basis for $\LL$. Then for $A\in\LL$, there is a unique point
$a=[a_1,a_2,\ldots,a_k]^\rt \in \R^k$ such that
\begin{equation}\label{eq:A linear}
   A=\sum_{i=1}^{k}a_i S_i.
\end{equation}
We write $A=g(a)$. Since $A$ is determined by $a$, we consider the
perturbation $\Delta a \in \R^k$ on $a$. Then we denote $\Delta
A=g(a+\Delta a)-g(a)=g(\Delta a)$, since $g$ is linear.

\begin{lemma}\label{F-expression: psi}
The Fr\'{e}chet derivative ${\bf D}\phi([a;b])$ of function $\phi$ defined
in~(\ref{phi}) is given by
\begin{equation}\label{eq:phi D}
{\bf D}{\phi}([a;b])=MP(A, \lambda)\left[\sfv_1, \sfv_2,\ldots, \sfv_k ,
  A^\rt\right],
\end{equation}
where $\sfv_i=-A^\rt S_{i} x_\lambda+ S^\rt_{i} r_\lambda$ for
$i=1,2,\ldots,k$.
\end{lemma}

\proof Let $\Delta A=g(\Delta a)$ and $\Delta b$ be perturbations on
$A=g(a)$ and $b$ respectively. Firstly, denoting
${\cal A} =(A+\Delta A)^\rt (A+\Delta A)+\lambda^2L^\rt L$
and recalling that $P(A, \lambda ) = (A^\rt A + \lambda^2 L^\rt L)^{-1}$,
we have
\begin{align*}
{\cal A}&=\left(A^\rt A+\lambda^2 L^\rt L\right)+
\left(A^\rt (\Delta A) +(\Delta A)^\rt A\right)+ (\Delta A)^\rt (\Delta A)\\
&=\left(A^\rt A+\lambda^2 L^\rt
L\right)\left[I_n+P(A,\lambda)\left(A^\rt (\Delta A) + (\Delta
A)^\rt A\right)+P(A,\lambda)\left( (\Delta A)^\rt (\Delta
A)\right)\right] .
\end{align*}
If $\|\Delta A\|$ is sufficiently small, then
$\left\|P(A,\lambda)\left(A^\rt (\Delta A)+ (\Delta A)^\rt A+
(\Delta A)^\rt (\Delta A)\right)\right\|<1$, from
Lemma~\ref{lemma:Banach}, $\cal A$ is nonsingular and its inverse
\begin{align}\label{eq: cal A}
{\cal A}^{-1}&=\left[I_n+P(A,\lambda)\left(A^\rt (\Delta A)+
(\Delta A)^\rt A\right)+P(A,\lambda)
\left((\Delta A)^\rt (\Delta A)\right)\right]^{-1}P(A,\lambda)\nonumber \\
&=P(A,\lambda)-P(A,\lambda)\left(A^\rt (\Delta A)+ (\Delta A)^\rt
A\right)P(A,\lambda)+\Oh(\|\Delta A\|^2),
\end{align}
since $\left\|A^\rt (\Delta A)+ (\Delta A)^\rt
A\right\|=\Oh(\|\Delta A\|)$ and $\left\|(\Delta A)^\rt (\Delta
A)\right\|=\Oh(\|\Delta A\|^2)$. From \eqref{deltax}, \eqref{eq: cal
A} and $x_\lambda=P(A,\lambda)A^\rt b$, after some algebraic
manipulation, we have
\begin{equation*}
\Delta x=P(A,\lambda)\Big(A^\rt (\Delta b)+ (\Delta
A)^\rt(b-Ax_{\lambda})-A^\rt (\Delta A)x_{\lambda}\Big)+\Oh(\|\Delta
A\|^2)+\Oh(\|\Delta A\|\|\Delta b\|).
\end{equation*}
Omitting the second and higher order terms and applying
the third equation in \eqref{eq:kr}, we have
\begin{align}
\Delta x&\approx P(A,\lambda)\Big(A^\rt (\Delta b)+ (\Delta
A)^\rt(b-Ax_{\lambda})-A^\rt (\Delta A)x_{\lambda}\Big)\nonumber\\
&=P(A,\lambda)\left[ \left(-(x_{\lambda}^\rt \otimes A^\rt)
+ \Big(r_\lambda^\rt \otimes I_n\Big)\Pi \right)
\vect(\Delta A)+A^\rt (\Delta b)\right]\nonumber\\
&=P(A,\lambda)\left[ \left(-(x_{\lambda}^\rt \otimes A^\rt)
+ (I_n
\otimes r_\lambda^\rt ) \right)\vect(\Delta A)+A^\rt (\Delta b)
\right] , \label{eq:delta x}
\end{align}
recalling that $r_{\lambda} = b - A x_{\lambda}$.

Since $\Delta A$ is a structured perturbation on $A$,
then $\Delta A=g(\Delta a)$, i.e., there exist
parameters $\Delta a_1,\Delta a_2,\ldots,\, \Delta a_k$ such that
$\Delta A =\sum\limits_{i=1}^{k}\Delta a_i S_i.$ Denote  $\Delta
a=\left[\Delta a_1,\Delta a_2,\ldots,\Delta a_k\right]^\rt$. From
\eqref{eq:delta x}, we have
\begin{align*}
& \phi([a+\Delta a; b+\Delta b])-\phi([a;b])\\
&\approx M P(A,\lambda)\left\{\left(-(x_{\lambda}^\rt \otimes A^\rt ) +
(I_n \otimes r_\lambda^\rt ) \right)
[\vect(S_{1}),\ldots,\vect(S_{k})]\Delta a+ A^\rt (\Delta b) \right\}\\
&=MP(A,\lambda)\left[-A^\rt S_{1} x_\lambda+  S^\rt_{1} r_\lambda, \ldots ,
-A^\rt S_{k} x_\lambda+ S^\rt_{k} r_\lambda, A^\rt\right]\Delta v,
\end{align*}
where $\Delta v=[\Delta a;\Delta b]$.
By the definition of the Fr\'{e}chet derivative,
the lemma then can be proved. \eproof

\begin{theorem}\label{theorem LS: linea struct}
Let $A \in \str$, $b \in \R^m$ and $x_\lambda=\left(A^\rt A+
\lambda^2 L^\rt L \right)^{-1} A^\rt b=P(A,\lambda)A^\rt b$ be the Tikhonov
regularized solution of (\ref{TikReg1}). Then we obtain the structured
normwise, componentwise, and mixed condition numbers:
\begin{eqnarray*}
m^{\reg}_\str (A,b)&=&\frac{\left\|\sum\limits_{i=1}^{k}|a_i|\left|
M P(A, \lambda)\left(A^\rt S_{i} x_\lambda- S^\rt_{i}
r_\lambda\right)\right|+\left|M P(A,\lambda)A^\rt
  \right||b|\right\|_\infty}{\|M x_\lambda\|_\infty},\\
c^{\reg}_\str(A,b)&=&\left\|\frac{\sum\limits_{i=1}^{k}|a_i|\left|M
P(A, \lambda)\left( A^\rt S_{i} x_\lambda-  S^\rt_{i}
r_\lambda\right)\right|+\left|M P(A,\lambda)A^\rt
\right||b|}{M x_\lambda}\right\|_\infty,\\
\kappa^{\reg}_\str(A,b)
&=&\frac{\left\| M P(A, \lambda)\left[S^\rt_{1} r_\lambda-
A^\rt S_{1} x_\lambda,\ldots , S^\rt_{k} r_\lambda-A^\rt S_{k} x_\lambda ,
\ A^\rt\right]\right\|_2
\left\|\left[\begin{matrix}a \cr b\end{matrix}\right]\right\|_2}
{\|M x_\lambda\|_2}.
\end{eqnarray*}
\end{theorem}

\proof From Lemmas~\ref{def:CN} and~\ref{F-expression: psi}, we have
\begin{align*}
m^{\reg}_\str(A,b)
&=\frac{\left\|~\left|{\bf D}{\phi}\left([a;b]\right)\right|\
\left[\begin{matrix}|a| \cr |b|\end{matrix}\right]\right\|_\infty}
{\|x_\lambda\|_\infty}\\
&=\frac{\left\|\left|MP(A,\lambda)\left[\sfv_1; \sfv_2;\ldots ;
\sfv_k\right]\right| |a|+
\left|MP(A,\lambda) A^\rt\right||b|\right\|_\infty}{\|Mx_\lambda\|_\infty} \\
&=\frac{\left\|\sum\limits_{i=1}^{k}
|a_i|\left| MP(A, \lambda)\left(A^\rt S_{i} x_\lambda-
S^\rt_{i} r_\lambda\right)\right|+
\left|MP(A,\lambda)A^\rt \right||b|\right\|_\infty}
{\|Mx_\lambda\|_\infty} .
\end{align*}
and

Similarly, we can obtain explicit expressions of the structured
componentwise and normwise condition numbers.
\eproof

When $\{ S_k \}$ is the canonical basis for $\R^{m \times n}$ in
Theorem~\ref{theorem LS: linea struct}, we have the following compact
forms of the unstructured condition numbers in Remark \ref{remark:1} for
$m^{\reg} (A,b),\, c^{\reg} (A,b)$ and $\kappa^{\reg} (A,b)$.
\begin{theorem}
As stated before, we have the following expressions
\begin{eqnarray*}
m^{\reg} (A,b)&=&\frac{\left\|\left|
M P(A, \lambda)\left[ (I_n\otimes r_{\lambda}^\rt ) -
(x^\rt \otimes A^\rt ) \right]\right|\vect(|A|)+\left|M P(A,\lambda)A^\rt
  \right||b|\right\|_\infty}{\|M x_\lambda\|_\infty},\\
c^{\reg}(A,b)&=&\left\|\frac{\left|
M P(A, \lambda)\left[ (I_n\otimes r_{\lambda}^\rt ) -
(x^\rt \otimes A^\rt ) \right]\right|\vect(|A|)+\left|M P(A,\lambda)A^\rt
\right||b|}{M x_\lambda}\right\|_\infty,\\
 \kappa^{\reg}(A,b)
  &=&\frac{\left\| M P(A, \lambda)\left[ (I_n\otimes r_{\lambda}^\rt ) -
  (x^\rt \otimes A^\rt ) ,
 A^\rt\right]\right\|_2\sqrt{\left\|A
 \right\|_F^2+\|b\|_2^2}}{\|M x_\lambda\|_2}.
\end{eqnarray*}
\end{theorem}

\proof For the expression of  $ \kappa^{\reg}_\str(A,b)$ given in
Theorem \ref{theorem LS: linea struct},
let $\{ S_{ij}={e_{i}}^{(m)}{e_{j}^{(n)}}^\rt \}$ be the canonical basis
for $\R^{m \times n}$, where ${e_{j}}^{(n)}$ is the $j$-th column
of the identity matrix $I_n$,
$i=1,2,\ldots, m$ and $j=1,2,\ldots, n$. Then we have the following
simplified expression:
$$
-  S^\rt_{ij}
r_\lambda+A^\rt
S_{ij} x_\lambda=-{e_{j}}^{(n)} r_{\lambda,(i)}+
A^\rt {e_{i}}^{(m)}x_{\lambda,(j)},
$$
where $r_{\lambda,(i)}$ and $x_{\lambda,(j)}$ are
respectively the $i$-th and $j$-th components of $r_\lambda$ and $x_\lambda$.
Now, fixing $j$, we get
$$
\left[-{e_{j}}^{(n)} r_{\lambda,(1)}+A^\rt {e_{1}}^{(m)}
x_{\lambda,(j)},\ldots, -{e_{j}}^{(n)} r_{\lambda,(n)}+
A^\rt {e_{n}}^{(m)}x_{\lambda,(j)}\right]=
-{e_{j}}^{(n)}r_{\lambda}^\rt +x_{\lambda,(j)} A^\rt,
$$
which implies that
\begin{align*}
&\left[-  S^\rt_{11}
r_\lambda+A^\rt
S_{11} x_\lambda,\ldots,-  S^\rt_{m1}
r_\lambda+A^\rt
S_{m1} x_\lambda, -  S^\rt_{12}
r_\lambda+A^\rt
S_{12} x_\lambda\ldots,-  S^\rt_{mn}
r_\lambda+A^\rt
S_{mn} x_\lambda\right]\\
&=\left[-({e_{1}}^{(n)}\otimes r_{\lambda}^\rt ) + x_{\lambda,(1)} A^\rt,
-({e_{2}}^{(n)}\otimes r_{\lambda}^\rt ) + x_{\lambda,(2)} A^\rt,\ldots,
-({e_{n}}^{(n)}\otimes r_{\lambda}^\rt ) + x_{\lambda,(n)} A^\rt  \right]\\
&=\left[- (I_n\otimes r_{\lambda}^\rt ) +
( x^\rt \otimes A^\rt ) \right].
\end{align*}
Applying the above equation to the expression of $\kappa_{\str}^{\reg}(A,b)$
in Theorem~\ref{theorem LS: linea struct},
we prove the third statement. The expressions of $m^{\reg} (A,b)$ and
$c^{\reg} (A,b)$ can be obtained similarly.
\eproof
\begin{remark}
If we choose $M=I_n$, then $m^{\reg} (A,b),\, c^{\reg} (A,b)$ and
$\kappa^{\reg} (A,b)$ respectively reduce to the expressions of
${\rm m_{Reg}},\, {\rm c_{Reg}}$ and ${\rm  cond^F_{Reg}}$ in
\eqref{mReg}, \eqref{cReg} and \eqref{condF}.
\end{remark}

How are the structured condition numbers compared to their
unstructured counterparts? The cases of nonsingular matrix inversion
and linear systems have been investigated in
\cite{BG,BF,Rump03a,Rump03b,Rump09} and the references therein.
In the following proposition, we will show that
$m^{\reg}_\str (A,b)$ is smaller than $m^{\reg} (A,b)$. The same is true
for the componentwise and normwsie condition numbers.
Before that we need the following lemma for rectangular structured matrices.
Its proof is omitted since it is similar to that of \cite[Lemma 6.3]{Rump03a}.

\begin{lemma}\label{lemma:a norm estimation}
when $A$ is a Toeplitz or Hankel matrix,
and $A=\sum\limits_{i=1}^{k}a_i S_i$, then
$$
\|a\|_2\leq \sqrt{2}\|A\|_F.
$$
\end{lemma}

\begin{proposition}\label{pro:relation}
Suppose that the basis $\{ S_1,S_2,\ldots, S_k \}$ for $\str$ satisfies
$|A|=\sum\limits_{i=1}^k |a_i||S_i|$ for any $A\in \str$
in \eqref{eq:A linear}, then
\[
m^{\reg}_\str (A,b)\leq m^{\reg} (A,b) \quad \hbox{and} \quad
c^{\reg}_\str (A,b)\leq c^{\reg} (A,b).
\]
For the structured normwise condition number, when $A$ is
a Toeplitz or Hankel matrix, we have
$$
\kappa^{\reg}_\str (A,b)\leq
\sqrt{2} \max\left\{\max_{i=1,2,\ldots,k}\|S_i\|_F,1\right\}
\kappa^{\reg} (A,b).
$$
\end{proposition}
\proof
From Theorem \ref{theorem LS: linea struct}, using the monotonicity of
the infinity norm, we have
\begin{align*}
& \left\|\sum\limits_{i=1}^{k}|a_i|
\left| M P(A, \lambda)\left(A^\rt S_{i} x_\lambda-
S^\rt_{i} r_\lambda\right)\right| +
\left|M P(A,\lambda)A^\rt \right||b|\right\|_\infty\\
&=\left\|\left|MP(A,\lambda)
\left[A^\rt S_{1} x_\lambda-  S^\rt_{1} r_\lambda, \ldots ,
A^\rt S_{k} x_\lambda- S^\rt_{k} r_\lambda\right] \right||a|+
\left|M P(A,\lambda)A^\rt \right||b| \right\|_\infty \\
&=\left\|\left|MP(A,\lambda)
\left[A^\rt S_{1} x_\lambda-  S^\rt_{1} r_\lambda, \ldots ,
A^\rt S_{k} x_\lambda- S^\rt_{k} r_\lambda, A^\rt\right]\right|
\begin{bmatrix}|a|\cr|b|\end{bmatrix}\right\|_\infty\\
&=\left\|\left|MP(A,\lambda)
\left[\left( ( x^\rt \otimes A^\rt ) -
( I_n\otimes r_\lambda^\rt ) \right)\vect(S_{1}),\ldots,
\vect(S_{k}), A^\rt\right]\right|
\begin{bmatrix}|a|\cr|b|\end{bmatrix}\right\|_\infty\\
&\leq\left\|\left[\left|MP(A,\lambda)
\left( ( x^\rt \otimes A^\rt ) -
( I_n\otimes r_\lambda^\rt ) \right)\vect(S_{1}),\ldots,
\vect(S_{k})\right|,\left|MP(A,\lambda) A^\rt\right|\right]
\begin{bmatrix}|a|\cr |b|\end{bmatrix}\right\|_\infty\\
&=\left\|\left|MP(A,\lambda)
\left[\left( ( x^\rt \otimes A^\rt ) -
( I_n\otimes r_\lambda^\rt ) \right)\vect(S_{1}),\ldots,
\vect(S_{k}), A^\rt\right]\right|
\begin{bmatrix}|a|\cr|b|\end{bmatrix}\right\|_\infty\\
&\leq\left\|\left|MP(A,\lambda)
\left( ( x^\rt \otimes A^\rt ) -
( I_n\otimes r_\lambda^\rt ) \right)\right|
\sum_{i=1}^k |a_i||\vect(S_{i})|+
\left|MP(A,\lambda) A^\rt\right||b|\right\|_\infty\\
&=\left\|\left| M P(A, \lambda)\left[ ( I_n\otimes r_{\lambda}^\rt ) -
( x^\rt \otimes A^\rt ) \right]\right|\vect(|A|)+
\left|M P(A,\lambda)A^\rt \right||b|\right\|_\infty,
\end{align*}
for the last equality we use the assumption $|A|=\sum_{i=1}^k |a_i||S_i|$.
With the above inequality, and the expressions of $m^{\reg}_\str (A,b)$,
$m^{\reg} (A,b)$, $c^{\reg}_\str (A,b)$, $c^{\reg} (A,b)$,
it is easy to prove the first two inequalities in this proposition.

When $A$ is a Toeplitz or Hankel matrix, the standard basis for
the Toeplitz matrix subspace or the Hankel matrix subspace is orthogonal
under the inner product $\langle B_1,B_2 \rangle=\tr \left(B_1^\rt
B_2\right)=[\vect(B_1)]^\rt \vect(B_2)$ for $B_1,\, B_2\in \R^{m\times n}$.
It is easy to deduce that
\begin{align*}
& \left\| M P(A, \lambda)
\left[S^\rt_{1} r_\lambda-A^\rt S_{1} x_\lambda,\ldots ,
S^\rt_{k} r_\lambda-A^\rt S_{k} x_\lambda, \  A^\rt\right]\right\|_2\\
&= \left\|MP(A,\lambda)
\left[\left( ( x^\rt \otimes A^\rt ) -
( I_n\otimes r_\lambda^\rt ) \right)\vect(S_{1}),\ldots, \vect(S_{k}),\
A^\rt\right]\right\|_2\\
&= \left\|MP(A,\lambda)
\left[\left( ( x^\rt \otimes A^\rt ) -
( I_n\otimes r_\lambda^\rt ) \right),\  A^\rt\right]
\begin{bmatrix}
\left[\vect(S_{1}),\ldots, \vect(S_{k})\right] & 0\cr
0 & I_m
\end{bmatrix} \right\|_2\\
&\leq \left\|MP(A,\lambda)
\left[\left( ( x^\rt \otimes A^\rt ) -
(I_n\otimes r_\lambda^\rt ) \right), \  A^\rt\right]
\right\|_2\left\|\begin{bmatrix}
\left[\vect(S_{1}),\ldots, \vect(S_{k})\right]& 0\cr
0 & I_m
\end{bmatrix}\right\|_2\\
&=\left\|MP(A,\lambda)
\left[\left( ( x^\rt \otimes A^\rt ) -
( I_n\otimes r_\lambda^\rt ) \right), \  A^\rt\right]\right\|_2
\max\left\{\max_{i=1,2,\ldots,k}\|S_i\|_F,1\right\},
\end{align*}
where for the last equation we used the orthogonality of the basis
$\{ S_i \}$. So from Lemma \ref{lemma:a norm estimation},
\begin{eqnarray*}
 & &\kappa^{\reg}_\str (A,b)
\leq\max\left\{\max_{i=1,2,\ldots,k}\|S_i\|_F, 1 \right\}
\frac{\left\|MP(A,\lambda)\left[\left( x^\rt \otimes A^\rt \right) -
\left( I_n\otimes r_\lambda^\rt \right), ~A^\rt\right]\right\|_2
\left\|\left[\begin{matrix}a
\cr b\end{matrix}\right]\right\|_2}{\|M x_\lambda\|_2} \\
&\leq& \max\left\{\max_{i=1,2,\ldots,k}\|S_i\|_F, 1 \right\}
\frac{\left\|MP(A,\lambda)\left[\left( x^\rt \otimes A^\rt \right) -
\left( I_n\otimes r_\lambda^\rt \right), ~A^\rt\right]\right\|_2
\sqrt{2\left\|A \right\|_F^2+\|b\|_2^2}}{\|M x_\lambda\|_2},
\end{eqnarray*}
which completes the proof of this proposition. \eproof

\begin{remark}
Clearly, the assumption $|A|=\sum_{i=1}^k |a_i||S_i|$ in
Proposition \ref{pro:relation} is satisfied for Toeplitz and
Hankel matrices.
\end{remark}

\subsection{Condition number estimators}

Efficiently estimating condition numbers is crucial in practice.
The condition number $\kappa^{\reg}_\str(A,b)$, for example,
involves the spectral norm of the $l\times (m+k)$ matrix
${\bf D}{\phi}([a;b])$, which can be expensive to compute
when $m$ or $k$ is large. The power
method can be used for fast condition number
estimation~\cite[page 289]{HighamBook}. Its major computation
is the matrix-vector multiplications ${\bf D}{\phi}([a;b]) h_1$ and
${\bf D}{\phi}([a;b])^\rt h_2$ with $h_1 \in \R^{m+k}$ and $h_2 \in \R^l$.
To consider ${\bf D}{\phi}([a;b])^\rt h$, for $h\in \R^l$ and
$\sfv_i$ defined in Lemma \ref{F-expression: psi}, denoting
$D = MP(A, \lambda )A^\rt$ and using
$P(A,\lambda)^\rt=P(A,\lambda)$, (\ref{eq:kr}) and
$\Pi^\rt=\Pi^{-1}$, we have
\begin{align}
\sfv_i^\rt P(A,\lambda)M^\rt h&= r_\lambda^\rt S_iP(A,\lambda)M^\rt
h-x_\lambda^\rt S_i^\rt A P(A, \lambda )M^\rt h\cr
&=\vect \left[r_\lambda^\rt
S_iP(A,\lambda)M^\rt h-x_\lambda^\rt S^\rt_iD^\rt h\right]\cr
&= \left[ \left(h^\rt MP(A,\lambda)\right) \otimes r_\lambda^\rt \right]
\vect(S_i) - \left[ (h^\rt D) \otimes x_\lambda^\rt \right]
\vect(S_i^\rt)\cr
&= \left[ \left(h^\rt MP(A,\lambda)\right) \otimes r_\lambda^\rt \right]
\vect(S_i)-\left[(h^\rt D)\otimes x_\lambda^\rt \right]\Pi
\vect(S_i)\cr
&=\vect(S_i)^\rt\left[(P(A,\lambda)M^\rt h) \otimes
r_\lambda-\Pi^\rt\left( D^\rt h\otimes x_\lambda \right)\right]\cr
&=\vect(S_i)^\rt\Pi^{-1}\left[\Pi\vect(r_\lambda h^\rt M
P(A,\lambda))-\vect(x_\lambda h^\rt D)\right]\cr &=[\Pi
\vect(S_i)]^\rt \vect\left[P(A,\lambda) M^\rt h r_\lambda^\rt  -
x_\lambda h^\rt D\right]\cr &=\left[\vect(S_i^\rt)^\rt\right]
\vect\left[P(A,\lambda) M^\rt h r_\lambda^\rt  - x_\lambda h^\rt
D\right]\cr &=\tr\left[S_i\left(P(A,\lambda) M^\rt h r_\lambda^\rt -
x_\lambda h^\rt D\right)\right]\label{eq:trace},
\end{align}
where we applied $[\vect(A_1)]^\rt \vect(A_2)=\tr \left(A_1^\rt
A_2\right)$ for the same dimensional matrices $A_1$ and $A_2$ in the
last equality. It follows from \eqref{eq:trace} that
\begin{equation}\label{eq:phi rt}
{\bf D}{\phi}([a;b])^\rt h =
\left(MP(A, \lambda)\left[\sfv_1, \ldots, \sfv_k ,
  A^\rt\right]\right)^\rt h=
\begin{bmatrix}\sfv_1^\rt P(A,\lambda)M^\rt h\cr
\vdots \cr \sfv_k^\rt  P(A,\lambda)M^\rt h\cr
D^\rt h\end{bmatrix}=
\begin{bmatrix}a(h)\cr D^\rt h\end{bmatrix},
\end{equation}
where $a(h)= [\tr \left( S_1( P(A,\lambda) M^\rt h r_\lambda^\rt -
x_\lambda h^\rt D ) \right), ...,
\tr \left( S_k ( P(A,\lambda) M^\rt h r_\lambda^\rt -
x_\lambda h^\rt D ) \right) ]^\rt$.
It leads to the following proposition.
\begin{proposition}\label{pro:adjoint}
The adjoint operator of ${\bf D}{\phi}([a;b])$, with the scalar products
$a_1^\rt a_2+b_1^\rt b_2$ and $h^\rt h$  in $\R^{k+m}$
and $\R^l$ respectively, is
\begin{equation}\label{eq:pro1}
{{\bf D} {\phi}([a;b])}^* :\  h \in \R^l \mapsto
\begin{bmatrix} a(h), & D^\rt h\end{bmatrix} \in \R^k\times\R^m .
\end{equation}
Furthermore, when $l=1$,
\begin{equation}\label{pro1}
\kappa^{\reg}_\str(A,b)=\frac{\sqrt{\sum\limits_{i=1}^k
s_i^2+\|D\|_2^2}\left\|\left[\begin{matrix}a
  \cr b\end{matrix}\right]\right\|_2}{\|M x_\lambda\|_2},
\end{equation}
where $s_i=\tr\Big(S_i\left(P(A,\lambda) M^\rt r_\lambda^\rt  -
x_\lambda  D\right)\Big), \, i=1,2,\ldots,k$.
\end{proposition}
\proof For any $(\Delta a,\Delta b)\in \R^{k}\times \R^m $ and $ h\in \R^l$,
from Lemma~\ref{F-expression: psi} and \eqref{eq:phi rt}, we have
\begin{eqnarray*}
& & \langle h,{\bf D}{\phi}([a;b])\cdot
(\Delta a ,\Delta b) \rangle = h^\rt\left( {\bf D}{\phi}([a;b])\cdot (\Delta a ,\Delta b)\right) = h^\rt {\bf D}{\phi}([a;b])
\begin{bmatrix}\Delta a\cr\Delta b\end{bmatrix}\\
&=& ({\bf D}{\phi}([a;b])^\rt h)^\rt
\begin{bmatrix}\Delta a\cr\Delta b\end{bmatrix}=a(h)^\rt (\Delta a) +(D^\rt h)^\rt (\Delta b)=\langle {\bf D}{\phi}([a;b])^* \cdot h,\
(\Delta a, \Delta b) \rangle,
\end{eqnarray*}
which proves the first part.
For the second part, noticing that
$$
\|{\bf D}{\phi}([a;b])\|_2=\|{\bf D}{\phi}([a;b])^\rt\|_2=\max_{h \neq
0}\frac{\left\|\left [a(h)^\rt, ~(D^\rt
h)^\rt \right]^\rt\right\|_2}{\|h\|_2}.
$$
and using (\ref{eq:phi rt}), where $h \in \R$ since $l=1$,
we can show that
$$
\|{\bf D}{\phi}([a;b])\|_2 =
\sqrt{\sum\limits_{i=1}^k s_i^2+\|D\|_2^2} ,
$$
which completes the proof. \eproof

\begin{remark}
When $l=1$, we compute
the conditioning of the $i$-th component of the solution.
In that case $M$ is the $i$-th
canonical vector of $\R^{1 \times n}$ and, in \eqref{pro1},
$P(A,\lambda)M^\rt$ is the $i$-th column of $P(A,\lambda)$ and
$D$ is the $i$-th row of $P(A,\lambda)A^\rt$.
\end{remark}

Using \eqref{eq:phi D} and \eqref{eq:phi rt}, we can now apply
the iteration of the power method~\cite[page 289]{HighamBook} in
Algorithm \ref{algorithm:1} to compute the normwise condition number
$\kappa^{\reg}_\str(A,b)$. In this algorithm, we assume that
$x_\lambda$, $r_\lambda$ and $\lambda$ are available.
When the GSVD (\ref{GSVD}) of $(A,L)$ is available,
a compact form of $P(A,\lambda)$ is given by
\begin{equation}\label{eq:p inver}
P(A,\lambda)=QR^{-1}
\begin{bmatrix}(\Sigma^2+\lambda^2 S^2)^{-1}& {\bf 0} \cr
{\bf 0} &I_{n-p}
\end{bmatrix}R^{-\rt}Q^\rt,
\end{equation}
which can be used to reduce the computational cost of the
estimators of the normwise, mixed and componentwise condition numbers.

\begin{algorithm}
\caption{The power method for estimating
$\kappa^{\reg}_\str(A,b)$ } \label{algorithm:1}
Select initial vector
$h\in \R^l$.
\begin{algorithmic}
\FOR {$p=1,2,\ldots $}
\STATE
Using \eqref{eq:p inver}, calculate $P(A,\lambda) M^\rt h
r_\lambda^\rt  - x_\lambda h^\rt D$. From \eqref{eq:phi rt},
denote $\bar a_p=a(h)$ and $\bar b_p=D^\rt h$. \\
\STATE
Calculate $\nu=\left\|[\bar a_p;~\bar b_p]\right\|_2$, let
$a_p=\bar a_p / \nu$ and $b_p=\bar b_p / \nu$.\\
\STATE
Let $A_p=\sum\limits_{i=1}^k a_{p,(i)} S_i$,
where $a_{p,(i)}$ is the $i$-th component of $a_p$. \\
\STATE Using \eqref{eq:phi
D} and \eqref{eq:p inver}, compute $h=MP(A,\lambda)\Big(A^\rt
b_p+ A_p^\rt r_\lambda -A^\rt A_p x_{\lambda}\Big)$.
 \ENDFOR \STATE $\kappa^{\reg}_\str(A,b)=\sqrt{\nu}$.
\end{algorithmic}
\end{algorithm}

The quantity $\nu$ computed by Algorithm \ref{algorithm:1} is
an approximation of the
largest eigenvalue of ${{\bf D} {\phi}([a;b])} {{\bf D} {\phi}([a;b])}^\top $.
When there is an estimate of the corresponding dominant eigenvector of
${{\bf D} {\phi}([a;b])} {{\bf D} {\phi}([a;b])}^\top $,
the initial $h$ can be set to this estimate, but in many implementations
$h$ is initialized as a random vector.
The algorithm is terminated by a sufficient number of iterations or
by evaluating the difference between two consecutive values of $\nu$
and comparing it to a tolerance
given by the user.

For the mixed and componentwise condition numbers,
we note that
\begin{align*}
m^{\reg}_\str(A,b)&=
\frac{\left\|MP(A, \lambda)\left[\sfv_1, \ldots, \sfv_k ,
A^\rt\right]\Diag([a;b])\right\|_\infty}{\|Mx_\lambda\|_\infty}=
\frac{\left\|{\bf D}{\phi}([a;b])\Diag([a;b])\right\|_\infty}
{\|Mx_\lambda\|_\infty},\\
c^{\reg}_\str(A,b)&=
\left\|\frac{MP(A, \lambda)\left[\sfv_1, \ldots, \sfv_k ,
A^\rt\right]\Diag([a;b])}{Mx_\lambda}\right\|_\infty=
\frac{\left\|{\bf D}{\phi}([a;b])\Diag([a;b])\right\|_\infty}
{\|Mx_\lambda\|_\infty}.
\end{align*}
The above equations show that we only need to estimate the infinity
norm of ${\bf D}{\phi}([a;b])\Diag([a;b])$. Since we have the adjoint
operator of ${\bf D}{\phi}([a;b])$ in \eqref{eq:pro1}, the power method
for estimating one norm~\cite[page 292]{HighamBook} can be used to
estimate $m^{\reg}_\str(A,b)$ as shownin Algorithm  2. Also, note
that, from \eqref{eq:phi rt}, for $h\in \R^l$,
\begin{equation} \label{eq:phiDiag}
\left({\bf D}{\phi}([a;b])\Diag([a;b])\right)^\rt h=
\Diag([a;b])\begin{bmatrix}
a(h)\cr D^\rt h\end{bmatrix}=
\begin{bmatrix}a\odot a(h)\cr b \odot (D^\rt h)\end{bmatrix},
\end{equation}
where `$\odot$' denotes the Hadamard (componentwise) product.
In Algorithm~2, ${\rm sign}(a)$ denotes the vector obtained
by applying the sign function to each component of the vector $a$.
We can estimate $c^{\reg}_\str(A,b)$ similarly.

\begin{algorithm}
\caption{The power method for estimating $m^{\reg}_\str(A,b)$ }
\label{al:mixed}
Select initial vector $h=l^{-1}\bfe\in \R^l$.
\begin{algorithmic}
\FOR {$p=1,2,\ldots $}
\STATE
Using \eqref{eq:p inver}, calculate
$P(A,\lambda) M^\rt h r_\lambda^\rt  - x_\lambda h^\rt D$.
From \eqref{eq:pro1}, compute $a(h)$ and $D^\rt h$. \\
\STATE
Using \eqref{eq:phiDiag}, denote $\alpha_p=a\odot a(h)$ and
$\beta_p=b\odot( D^\rt h)$.
\STATE
Let $\bar a_p={\rm sign}\left(\alpha_p\right)$ and
$\bar b_p={\rm sign}\left(\beta_p\right)$.
\STATE
Compute $ a_p=a\odot {\bar a_p},\, b_p=b\odot {\bar b_p}$.
\STATE
Form $A_p=\sum\limits_{i=1}^k a_{p,(i)} S_i$,
where $a_p=[a_{p,(1)},a_{p,(2)},\ldots,a_{p,(k)}]^\rt$.
\STATE
Using \eqref{eq:phi D} and \eqref{eq:p inver}, compute
$z=MP(A,\lambda)(A^\rt b_p+ A_p^\rt r_\lambda -A^\rt A_p x_{\lambda})$.
\IF
{ $\|z\|_\infty \leq h^\rt z$}
\STATE
$\gamma=\left\|\begin{bmatrix}
\alpha_p \cr \beta_p
\end{bmatrix}\right\|_1$
\STATE
quit
\ENDIF
\STATE
$h=e_j^{(l)}$, where $|z_j|=\|z\|_\infty$ (smallest such $j$).
\ENDFOR
\STATE
$m^{\reg}_\str(A,b)=\gamma$.
\end{algorithmic}
\end{algorithm}

The main computational cost of Algorithm~\ref{algorithm:1}
or Algorithm~\ref{al:mixed} is the computation
of solving several nonsingular triangular systems with the coefficient
matrices $R$ and $R^{\rt}$. If we have the GSVD of $(A,L)$ available, the
computational cost is insignificant compared with the cost of solving the
Tikhonov regularized problem. Thus, the estimators can be integrated
into a GSVD based Tikhonov solver without compromising the
overall computational complexity. Our methods can be readily modified for
fast unstructured condition number estimation, which is not
considered in \cite{CLTW}.

\section{Nonlinear Structures}\label{sec:nonlinear}

In this section, we present the structured condition numbers
of matrices with nonlinear structures, namely the Vandermonde
matrices and the Cauchy matrices.

\subsection{Vandermonde matrices}

Let $\VdM$ be the class of $m\times n$ Vandermonde matrices. If
$V=[v_{ij}]\in\VdM$, then there exists
$a=[a_0,a_1,\ldots,a_{n-1}]^\rt\in\R^n$ such that, for all
$i=0,1,\ldots,m-1$ and $j=0,1,\ldots,n-1$, $v_{ij}=a_j^i$. We write
$V=g(a)$. Let $\Delta a=(\Delta a_0,\Delta a_1,\ldots,\Delta
a_{n-1})^\rt\in\R^n$  be the perturbation on $a$. Then we define the
first order term $\Delta V$ of $g(a+\Delta a)-g(a)$.  From

\begin{lemma} \label{Express delta: Vandermonde}
{\rm \bf (\cite[Lemma 6]{CD07})}
An explicit expression of $\Delta V$ is
$$
\Delta V=V_1\Diag(\Delta a), \quad \hbox{where} \quad
V_1=\Diag(c)\left[\begin{matrix}{\bf  0}\cr
V(1:m-1,:)\end{matrix}\right], \ c=[0,1,\ldots,m-1]^\rt .
$$
Here $V(1:m-1,:)$ is the $(m-1) \times n$  submatrix of $V$
consisting of the first $m-1$ rows of $V$.
\end{lemma}

\begin{lemma}\label{F-expression: Vdm}
The Fr\'{e}chet derivative ${\bf D} \phi ([a;b])$ of
function $\phi$ defined in~(\ref{phi}) is
$$
{\bf D} {\phi}([a;b])=M P(V, \lambda) \left[-V^\rt
V_1\Diag(x_\lambda)+\Diag(y), V^\rt\right],
$$
\end{lemma}
where $y=V_1^\rt r_\lambda$ and $r_\lambda=b-Vx_\lambda$.

\proof It follows
from \eqref{eq:delta x}  and $\Diag(a)z=\Diag(z)a$ for
vectors $a$ and $z$ of the same dimension,
\begin{align*}
&\phi([a+\Delta a; b+\Delta b])-\phi([a;b])\\
&\approx MP(V,\lambda)\left[ \left(-x_{\lambda}^\rt \otimes V^\rt+
I_n \otimes r_\lambda^\rt \right)\vect(\Delta V)+V^\rt (\Delta b)\right]\\
&=MP(V,\lambda)\left[ \left(-x_{\lambda}^\rt \otimes V^\rt+
(r_\lambda^\rt\otimes I_n)\Pi \right)\vect(\Delta V)+
V^\rt (\Delta b)\right]\\
&=MP(V,\lambda)\left[ \left(-  V^\rt (\Delta V)x_{\lambda}+
(\Delta V)^\rt r_\lambda \right)+V^\rt (\Delta b)\right]\\
&=M P(V, \lambda)\left[\left(-V^\rt V_1\Diag(\Delta a) x_\lambda+
\Diag(\Delta a)V_1^\rt r_\lambda\right)+V^\rt (\Delta b)\right]\\
&= M P(V, \lambda) \left[-V^\rt V_1\Diag(x_\lambda)+
\Diag(y),\  V^\rt\right]
\left[\begin{matrix}\Delta a \cr \Delta b
\end{matrix}\right],
\end{align*}
which completes the proof of this lemma. \eproof

From Lemmas~~\ref{def:CN} and \ref{F-expression: Vdm}, we
have the following theorem of structured condition numbers
of the Vandermonde matrix.

\begin{theorem}\label{theorem VdM: linea struct}
Let $V \in \VdM$, $b \in \R^m$ and $x_\lambda=\left(V^\rt V+\lambda
L^\rt L\right)^{-1} V^\rt b=P(V,\lambda)V^\rt b$ be the Tikhonov
regularized solution of (\ref{TikReg1}). Recall that $y=V_1^\rt
r_\lambda$, then the structured condition numbers of the
Vandermonde matrix are:

\begin{eqnarray*}
m^{\reg}_\VdM (V,b)&=&\frac{\left\|\left|M P(V, \lambda)
\left[V^\rt V_1\Diag(x_\lambda)-\Diag(y)\right]\right| |a|+
\left|P(V,\lambda)V^\rt\right||b|\right\|_\infty}{\|M x_\lambda\|_\infty},\\
c^{\reg}_\VdM(V,b)&=&\left\|\frac{\left|M P(V, \lambda)
\left[V^\rt V_1\Diag(x_\lambda)-\Diag(y)\right]\right| |a|+
\left|P(V,\lambda)V^\rt\right||b|}{M x_\lambda}\right\|_\infty,\\
\kappa^{\reg}_\VdM(V,b) &=&\frac{\left\|M P(V, \lambda) \left[\Diag(y)-
V^\rt V_1\Diag(x_\lambda), V^\rt\right]\right\|_2
\left\|\left[\begin{matrix}
a \cr b\end{matrix}\right]\right\|_2}{\|M x_\lambda\|_2}.
\end{eqnarray*}
\end{theorem}
In particular, when $l=1$,
\[
\kappa^{\reg}_\VdM(V,b)=
\frac{\sqrt{\left\|y \odot \left( P(V,\lambda)M^\rt \right) -
x_\lambda \odot \left( V_1^\rt D_V^\rt \right) \right\|_2^2+
\|D_V\|_2^2} \left\| \left[ \begin{matrix}
a \cr b\end{matrix}\right]\right\|_2}{\|M x_\lambda\|_2}.
\]

Analogous to Proposition~\ref{pro:adjoint},
the adjoint operator of ${\bf D}{\phi}([a;b])$, using the scalar products
$a_1^\rt a_2+b_1^\rt b_2$ and $h^\rt h$ on $\R^{m+n}$
and $\R^l$ respectively, is
\[
{\bf D}{\phi}([a;b])^* :\ h \in \R^l \mapsto
[y \odot P(V,\lambda)M^\rt h)-
x_\lambda \odot(V_1^\rt D_V^\rt h) \  D_V^\rt h] \in \R^n\times\R^m ,
\]
where $D_V=MP(V,\lambda)V^\rt$.
The above expressions can be used to to estimate $m^{\reg}_\VdM (V,b)$,
$c^{\reg}_\VdM (V,b)$ and $\kappa^{\reg}_\VdM (V,b)$ with lower
dimensional input. We can devise algorithms similar to
Algorithms~\ref{algorithm:1} and \ref{al:mixed} for estimating
the condition numbers.

\subsection{Cauchy matrices}

Let $\Cauchy$ be the class of $m\times n$ Cauchy matrices. If
$C=[c_{ij}] \in\Cauchy$, then there exist $u=[u_1,u_2,\ldots,
u_m]^\rt\in \R^m$ and $v=[v_1,v_2,\ldots, v_n]^\rt\in \R^n$, with
$u_i \neq v_j$ for $i=1,2,\ldots,m$, $j=1,2,\ldots,n$ such that, for
all $i\leq m$ and $j\leq n$,
$$
c_{ij}=\frac{1}{u_i-v_j}.
$$
If $w=[u;v]\in \R^{m+n}$, then $C=g(w)$. Let $\Delta w=[\Delta u;
\Delta v]=[\Delta u_1,\ldots,\Delta u_m,\Delta v_1, \ldots,
\Delta v_n]^\rt$ $\in\R^{m+n}$ be the perturbation on
$w$. The first order term
$\Delta C$ in $g(w+\Delta w)-g(w)$ is given by \cite[Lemma 9]{CD07}
$$
\Delta C \approx \left[\frac{\Delta u_i -\Delta v_j}{(u_i-v_j)^2}\right]
=\Diag(\Delta u)C_1-C_1\Diag(\Delta v) \in \R^{m \times n},
$$
where $C_1= [ 1 / (u_i-v_j)^2 ] \in \R^{m \times n}$.

\begin{lemma}\label{lemma: Cauchy de}
The Fr\'{e}chet derivative ${\bf D}\phi([w;b])$  of function $\phi$
defined in~(\ref{phi}) is given by
$$
{\bf D}{\phi}([w;b])=M P(C, \lambda)\left[C_u,C_v, C^\rt\right],
$$
where $C_u=C_1^\rt \Diag( r_\lambda)-C^\rt \Diag(z_1),~ C_v= C^\rt
C_1\Diag(x_\lambda)-\Diag(z_2),~z_1=C_1 x_\lambda$, $z_2=C_1^\rt
r_\lambda$ and $r_\lambda=b-Cx_\lambda$. \end{lemma}

\proof
Following the proof of Lemma~\ref{F-expression: Vdm}, we can show that
\begin{eqnarray*}
& & \phi([w+\Delta w; b+\Delta b])-\phi([w;b])\\
&=& M P(C, \lambda)\left[C_1^\rt \Diag( r_\lambda)-
C^\rt \Diag(z_1),\  C^\rt C_1 \Diag(x_\lambda)- \Diag(z_2),\
C^\rt \right]
\left[\begin{matrix}\Delta w \cr \Delta b
\end{matrix}\right].
\end{eqnarray*}

Then the Fr\'{e}chet derivative of $\phi$ at $[w;b]$ is
$$
{\bf D}{\phi}[w;b]=M P(C, \lambda)\left[C_1^\rt \Diag( r_\lambda)-C^\rt
\Diag(z_1), C^\rt C_1\Diag(x_\lambda)-\Diag(z_2), C^\rt\right].
$$

\begin{theorem}\label{theorem Cauchy: linea struct}
Let $C \in \Cauchy$, $b \in \R^m$ and
$x_\lambda=\left(C^\rt C+\lambda L^\rt L\right)^{-1} C^\rt b=
P(C,\lambda)C^\rt b$ be the solution of the
Tikhonov regularization problem (\ref{TikReg1}), then the structured
condition numbers are:
\begin{eqnarray*}
m^{\reg}_\Cauchy (C,b)&=&
\frac{\left\|\left|M P(C,\lambda)C_u\right| |u|+
\left|M P(C,\lambda)C_v\right| |v|+
\left|M P(C,\lambda)C^\rt\right||b|\right\|_\infty}
{\|M x_\lambda\|_\infty},\\
c^{\reg}_\Cauchy(C,b)&=&
\left\|\frac{\left|M P(C,\lambda)C_u\right| |u|+
\left|M P(C,\lambda)C_v\right| |v|+
\left|M P(C,\lambda)C^\rt\right||b|}{M x_\lambda}\right\|_\infty,\\
\kappa^{\reg}_\Cauchy(C,b) &=&
\frac{\left\|M P(C, \lambda)\left[C_u,C_v, C^\rt\right]\right\|_2
\left\| \left[\begin{matrix}w \cr b
\end{matrix}\right]\right\|_2}{\|M x_\lambda\|_2}.
\end{eqnarray*}
\end{theorem}
In particular, when $l=1$,
\[
\kappa^{\reg}_\Cauchy(C,b)=
\frac{\sqrt{t^2+s^2+\|D_C\|_2^2}}{\|M x_\lambda\|_2}
\left\|\left[\begin{matrix}u\cr v \cr b
\end{matrix}\right]\right\|_2,
\]
where $t=\left\|r_\lambda \odot\left[C_1P(C,\lambda)M^\rt \right]-
z_1\odot (D_C^\rt ) \right\|_2$ and
$s=\left\|x_\lambda \odot (C_1^\rt D_C^\rt )-
z_2\odot \left[P(C,\lambda)M^\rt \right]\right\|_2$.

Similar to the case of the Vandermonde matrix, for the Cauchy matrix,
the adjoint operator of ${\bf D}{\phi}([w;b])$, using the scalar products
$u_1^\rt u_2+v_1^\rt v_2+b_1^\rt b_2$ and $h^\rt h$ on
$\R^{2m+n}$ and $\R^l$ respectively, is
\[
{\bf D}{\phi}([w;b])^* : \  h \in \R^l \mapsto
\left[ u(h)\  v(h) \  D_V^\rt h \right] \in \R^m\times \R^n\times\R^m ,
\]
where $u(h)=r_\lambda \odot \left[C_1P(C,\lambda)M^\rt h\right]-
z_1\odot (D_C^\rt h)$, $v(h)=x_\lambda \odot(C_1^\rt D_C^\rt h)-
z_2\odot \left[P(C,\lambda)M^\rt h\right]$ and
$D_C=MP(C,\lambda)C^\rt$.

In particular, when $l=1$, we have
\begin{equation*}
\kappa^{\reg}_\Cauchy(V,b)=\frac{\sqrt{t^2+s^2+\|D_C\|_2^2}}{\|M
x_\lambda\|_2}\left\|\left[\begin{matrix}u\cr v
  \cr b\end{matrix}\right]\right\|_2,
\end{equation*}
where $t=\left\|r_\lambda \odot\left[C_1P(C,\lambda)M^\rt
\right]-z_1\odot (D_C^\rt ) \right\|_2$ and $s=\left\|x_\lambda
\odot(C_1^\rt D_C^\rt )-z_2\odot \left[P(C,\lambda)M^\rt
\right]\right\|_2$.

Using the above expressions,
the algorithms similar to  Algorithms~\ref{algorithm:1} and \ref{al:mixed}
for estimating $m^{\reg}_\Cauchy (C,b)$, $c^{\reg}_\Cauchy (C,b)$ and
$\kappa^{\reg}_\Cauchy (C,b)$ can be obtained.

\section{SCE for the Tikhonov Regularization Problem}\label{sec:SCE for Tik}

In this section we use SCE to devise algorithms for the condition
estimations of the structured and unstructured Tikhonov
regularization problem, both the normwise and componentwise cases
are considered.

\subsection{SCE for normwise perturbations}

For the unstructured Tikhonov regularization problem,
we are interested in the condition estimation for the function
$\psi([A,b])$ at the point $[A, b]$ defined in \eqref{psi}.
Let $[A\, \,b]$ be perturbed to $[A+\delta E\, \, b+\delta f]$ in
the normal equations (\ref{NormEqn}), where $\delta \in \R$,
$E\in \R^{m\times n}$ and $f \in \R^m$ and $[E\,\, f]$ has the
Frobenius norm equal to one. According to Subsection \ref{sec:SCE},
we first need to evaluate the directional derivative
${\bf D}\psi([A,b];[E,f])$ of $\psi([A,b])$ with respect to $[A,b]$
in the direction $[E,f]$.
From the proof of Lemma~\ref{F-expression: psi}, we have
$$
{\bf D}\psi([A,b];[E,f])=
P(A,\lambda)\Big(A^\rt f+ E^\rt r_{\lambda}-A^\rt Ex_{\lambda}\Big).
$$
When we have the GSVD (\ref{GSVD}) of $(A,L)$,  it is easy to deduce that
\begin{equation}\label{eq:directional derivative}
{\bf D}\psi([A,b];[E,f])=
QR^{-1}\begin{bmatrix}(\Sigma^2+\lambda^2 S^2)^{-1}& {\bf 0} \cr
{\bf 0} &I_{n-p}\end{bmatrix}R^{-\rt}Q^\rt\Big(A^\rt f+
E^\rt r_{\lambda}-A^\rt Ex_{\lambda}\Big).
\end{equation}
With the above result, we now use the results of Subsection~\ref{sec:SCE} to
obtain the SCE-based methods for estimating the
condition of the Tikhonov regularization problems.
Both the normwise and componentwise perturbations are considered.
Algorithm~\ref{algo:subnorm} computes an estimation of the
normwise condition number. Inputs
to the method are the matrices $A \in  \Rmn$, $L\in \Rpn$,
the vector $b \in \Rm$, the computed solution $x_\lambda$ and
the parameter $\lambda$. The output is
an estimation $\kappa_{\rm SCE}^{(k)}$ of the normwise condition number
${\rm  cond^F_{Reg}}$.
The method requires the GSVD (\ref{GSVD})
of $(A, L)$, which is generally computed when solving
the Tikhonov regularization problem. The integer $k\geq 1$ refers to
the number of perturbations of input data. Note that
when $k = 1$, there is no need to orthonormalize the set of
vectors in Step 1 of the
method. In the following the standard normal distribution are denote
by ${\cal N}(0, 1)$, and for $B=(b_{ij})\in \Rpq$,
$|B|^2=(|b_{ij}|^2)\in \Rpq$ and $\sqrt{|B|}=(\sqrt{|b_{ij}|})\in \Rpq$.

\begin{algorithm}%
\caption{SCE for the Tikhonov regularization problem
under normwise perturbations}
\label{algo:subnorm}
\begin{itemize}
\item[1.] Generate matrices $[E_1,~f_1],
[E_2,~f_2],\ldots, [E_k,~f_k]$ whose entries are random numbers in
${ \cal N}(0,1)$, where $E_i\in \R^{m\times n},\, f_i \in \Rm$. Use a QR
factorization for the matrix
$$
\left[\begin{matrix}\vect(E_1)&\vect(E_2)&\cdots&
\vect(E_k)\cr f_1&f_2&\cdots&
f_k
\end{matrix}\right]
 $$
and form an orthonormal matrix $[q_1,q_2,\ldots,q_k]$. Each $q_i$
can be converted into the desired matrices $[\widetilde{E_i}$,
$\widetilde{f_i} ]$ with the {\sf unvec} operation.
\item[2.] Calculate ${\bf D}\psi([A,b];[\widetilde{E_i},\widetilde{f_i}])$ by
\eqref{eq:directional derivative}, $i=1,2,\ldots,k$.
\item[3.] Compute the absolute condition vector
\begin{eqnarray*}
\kappa_{\rm
abs}^{(k)}&:=&\frac{\omega_k}{\omega_p}
\sqrt{\left|{\bf D}\psi([A,b];[\widetilde{E_1},\widetilde{f_1}])\right|^2+
\cdots+\left|{\bf D}\psi([A,b];[\widetilde{E_k},\widetilde{f_k}])\right|^2}.
\end{eqnarray*}
\item[4.] Compute the normwise condition estimation:
$$
\kappa_{\rm SCE}^{(k)}:=\frac{\left\|\kappa_{\rm
abs}^{(k)}\right\|_2 \left\|[A,~b]\right\|_F}{||x_\lambda||_2}.
$$
\end{itemize}
\end{algorithm}

\subsection{SCE for componentwise perturbations}

Componentwise perturbations are relative to the magnitudes of
the corresponding entries in the input arguments
(e.g., the perturbation $\Delta A$ satisfies $|\Delta A| \leq \epsilon |A|$,
see (\ref{mReg})). These perturbations
may arise from input error or from rounding error, and hence are
the most common perturbations encountered in practice. In fact,
most of error bounds in LAPACK are componentwise since
the perturbations of input
data are componentwise in real world computing, see
\cite[section 4.3.2]{Bai} for details.
We often want to find the condition of a function with respect to
componentwise perturbations on inputs. For the function
$$
\psi([A,~b])=\left(A^\rt A+\lambda^2 L^\rt L\right)^{-1}A^\rt b,
$$
SCE is flexible enough to accurately gauge the sensitivity of
matrix functions subject to componentwise perturbations.
Define the linear function
$$
h([B,~d])=[B,~d]\odot [A,~b],\quad B\in \R^{m\times n},\quad d\in \Rm.
$$
Let ${\cal E} \in \R^{m\times (n+1)}$ be
the matrix of all ones, then $h({\cal E})=[A~b]$ and
$$
h({\cal E}+[E,~f])=[A,~b]+h([E,~f]).
$$
We know that $h([E,~f])$ is a componentwise perturbation on $[A,~ b]$,
and $h$ converts a general perturbation $\cal E$ into
componentwise perturbations on $[A,~b]$. Therefore, to obtain
the sensitivity of the solution with respect
to relative perturbations, we simply evaluate the Fr\'{e}chet derivative of
$$
\psi([A,~b])=\psi(h({\cal E}))
$$
with respect to $\cal E$ in the direction $[E,~ f]$, which is
$$
{\bf D}(\psi \circ h)\left({\cal E};[E,~f])\right)=
{\bf D}\psi (h({\cal E})){\bf D}h\left({\cal E};[E,~f]\right)=
{\bf D}\psi ([A,~b])h\left([E,~f]\right)=
{\bf D}\psi \left([A,~b];h\left([E,~f]\right)\right),
$$
since $h$ is linear. Thus, to estimate the condition of the Tikhonov regularization solution $x_\lambda$
when perturbations are componentwise, we first generate the perturbations $E$
and $f$ and multiply them componentwise by the entries of $A$ and $b$,
respectively. The remaining steps are the same as
the corresponding steps in Algorithm~\ref{algo:subnorm},
as shown in Algorithm~\ref{algo:subcomp}.

\subsection{SCE for structured perturbations}\label{sec:sub structure}

The SCE also is flexible for
the condition estimation for structured Tikhonov regularization problem.
We are interested in the condition estimation for the function
$\phi$ defined in
\eqref{phi}, which defines the general function for structured
Tikhonov regularization problem. Because SCE can estimate the condition
of the each component of $x_\lambda$, we only need to choose $M=I_n$ in
\eqref{phi}.

\begin{algorithm}[t]
\caption{SCE for the Tikhonov regularization problem
under componentwise perturbations}
\label{algo:subcomp}
\begin{itemize}
\item[1.]
Generate matrices $[E_1,~f_1],
[E_2,~f_2],\ldots, [E_k,~f_k]$ whose entries are random numbers
in ${ \cal N}(0,1)$, where $E_i\in \R^{m\times n},\, f_i \in \Rm$.
Use a QR factorization for the matrix
$$
\left[\begin{matrix}\vect(E_1)&\vect(E_2)&\cdots&
\vect(E_k)\cr f_1&f_2&\cdots&
f_k
\end{matrix}\right]
$$
to form an orthonormal matrix $[q_1,q_2,\ldots,q_k]$. Each $q_i$
can be converted into the desired matrices $[E_i~f_i]$ with the {\sf unvec}
operation.
\item [2.]
For $i =1,\ldots, k$, set $[\widetilde{E_i},~\widetilde{f_i} ]$
to the componentwise product of $[A,~ b]$ and $[E_i,~f_i]$.
\item[3.]
Calculate ${\bf D}\psi([A,b];[\widetilde{E_i},\widetilde{f_i}])$ by
\eqref{eq:directional derivative}, $i=1,2,\ldots,k$.
\item[4.]
Compute the absolute condition vector
\[
c_{\rm abs}^{(k)} :=
\frac{\omega_k}{\omega_p}
\sqrt{\left|{\bf D}\psi([A,b];[\widetilde{E_1},\widetilde{f_1}])\right|^2+
\cdots+\left|{\bf D}\psi([A,b];[\widetilde{E_k},\widetilde{f_k}])\right|^2}.
\]
\item[5.]
The mixed condition estimation
$m_{\rm SCE}^{(k)}$ and componentwise condition estimation
$c_{\rm SCE}^{(k)}$ are:
\[
m_{\rm SCE}^{(k)} :=
\frac{\left\|c_{\rm abs}^{(k)}\right\|_{\infty}}{\|x_\lambda\|_{\infty}}
\quad \hbox{and} \quad
c_{\rm SCE}^{(k)} :=
\left\|\frac{c_{\rm
abs}^{(k)}}{x_\lambda}\right\|_{\infty}.
\]
\end{itemize}
\end{algorithm}

The key step in the SCE is the computation of
the directional derivative ${\bf D}\phi([a;b];[e;f])$ of
$\phi([a;b])$ with respect to $[a;b]$
in the direction $[e;f]$, where $e\in \R^{k}$ and $f \in \R^m$.
We have derived the explicit expressions of the Fr\'{e}chet
derivative ${\bf D}\phi([a;b])$ in Lemmas~\ref{F-expression: psi},
\ref{F-expression: Vdm} and \ref{lemma: Cauchy de} for a general linear
structure, Vandermonde or Cauchy matrix. Based on
Lemmas~\ref{F-expression: psi}, \ref{F-expression: Vdm} and
\ref{lemma: Cauchy de}, the three directional derivatives
${\bf D}\phi([a;b];[e;f])$ are:
$$
{\bf D}\phi([a;b];[e;f]) = P(A,\lambda)
\left(A^\rt f+ E^\rt r_{\lambda}-A^\rt E x_{\lambda}\right), \quad
e=(e_i) \in \R^k,\  E=\sum_{i=1}^k e_iS_i,
$$
for linear structures,
$$
{\bf D}\phi([a;b];[e;f]) = M P(V, \lambda) \left(\Diag(y) e-V^\rt
V_1\Diag(x_\lambda)e+ V^\rt f\right), \quad e\in \R^n,
$$
for Vandermonde matrices, and
$$
{\bf D}\phi([a;b];[e;f]) = M P(C, \lambda )
\left( C_u e_1+C_v e_2+C^\rt f \right), \quad
e= \left[\begin{matrix}e_1 \cr  e_2 \end{matrix}\right] \in \R^{m+n},
$$
for Cauchy matrices,
where $V_1$ is defined in Lemma~\ref{Express delta: Vandermonde},
$y=V_1^\rt r_\lambda$, $C_u$ and $C_v$ are defined in
Lemma~\ref{lemma: Cauchy de}. Based on those expressions,
we can derive algorithms for structured
normwise and componentwise condition estimation. The algorithms
are similar to those of Algorithms \ref{algo:subnorm} and
\ref{algo:subcomp}, thus are omitted here.

\section{Numerical Examples}\label{sec:numerical exam}

In this section, we demonstrate our test results of some
numerical examples to illustrate
structured condition numbers and condition estimations presented in
the previous sections. All the
computations are carried out using \textsc{Matlab} 8.1 with the
REGULARIZATION TOOLS package \cite{tool} with the machine precision
$2.2 \times 10^{-16}$.

For a structured matrix $A$, which is determined by the
vector $a\in \R^k$, we generated the perturbed matrix $\widehat A$ as follows.
For $a\in \R^k$ and $ b \in \R^m$, let
$[s;f]$ be a random vector whose entries are
uniformly distributed in
the open interval $(-1,1)$, where $s\in \R^k$ and $f\in \R^m$,
the perturbations on $a$ and $b$ are respectively
\begin{equation}\label{eq:sec numerical example perturbation}
\Delta a_{i}=\varepsilon s_{i}a_{i},  \qquad \Delta
b_{j}=\varepsilon f_{j}b_{j},
\end{equation}
then $\widehat A=g(a+\Delta a)$ and $\widehat b=b+\Delta b$.
In our experiments, we set $\varepsilon=10^{-8}$.

REGULARIZATION TOOLS package \cite{tool} includes four methods for
determining the Tikhonov regularization parameter. For the Tikhonov
regularization with continuous regularization parameter,
the $L$-curve is a continuous curve as a parametric plot of the
discrete smoothing (semi) norm $\|Lx_\lambda\|_2$ versus the
corresponding residual norm $\|Ax_\lambda-b\|_2$, with the parameter
$\lambda$ as the parameter. The corner of the $L$-curve
appears for regularization parameters close to the optimal parameter
that balances the regularization errors and perturbation errors in
$x_\lambda$, which is the basis for the $L$-curve criterion for
choosing the regularization parameter. Besides the $L$-curve criterion
for parameter-choice, a variety of parameter-choice strategies have
been proposed, such as the \textit{discrepancy principle} (Discrep. pr.)
\cite{morozov}, \textit{generalized cross-validation} (GCV)
\cite{wahba} and the quasi-optimality criterion (Quasi-opt) \cite{morozov}.

The Tikhonov regularization solution $x_\lambda$ was computed by the
Matlab function
{\sf tikhonov} corresponding to $A,\, b$ in REGULARIZATION TOOLS
package with different regularization parameters chosen by
four classical criteria or by the predefined value. The perturbed
solution $y_\lambda$ is obtained in
the similar way to $x_\lambda$, but $y_\lambda$ corresponds to
$\widehat A$ and $\widehat b$. Denote the error $\Delta
x_\lambda=y_\lambda-x_\lambda$.

We compare the structured condition numbers with unstructured ones
for various Tikhonov regularization parameters in the following
examples.

\begin{example}~({\bf \cite{Rump03b}})\label{ex:vande}
Let $A=g(a)$ be a $5\times 5$ symmetric Toeplitz matrix which is
defined 
$$
  A=g(a)=\begin{bmatrix}
    0&0&1+h&-1&1\cr
    0&0&0&1+h&-1\cr
    1+h&0&0&0&1+h\cr
    -1&1+h&0&0&0\cr
    1&-1&1+h&0&0
  \end{bmatrix},\quad a=\begin{bmatrix}0\cr 0\\ 1+h\\-1\\1\end{bmatrix}, \quad
  b=\begin{bmatrix}0\\h\\2(1+h)\\h\\ 0\end{bmatrix}
$$
for $h=10^{-3}$.

The above matrix $A$ is a square symmetric Toeplitz matrix
and $g(a)$ is a square symmetric Toeplitz matrix whose first
column is $a$. We can choose the
basis $Z_i=g( e_i )$, $i=1,2,\ldots,5$, so that $T=\sum_{i=1}^{5}a_i Z_i$.
From Theorem \ref{theorem LS: linea struct}, with the $Z_i$,
we can get the expressions for
$m_{\sf SymToep}(A,b)$, $c_{\sf SymToep}(A,b)$ and
$\kappa_{\sf SymToep}(A,b)$. The relative errors and condition numbers
are shown in Table \ref{Ex:toep}, where $M=I_5$.

\begin{table}\centering
\caption{\label{Ex:toep}  $L=I_5$, $M=I_5$.}
\begin{tabular}{||c|c c c c||}
\hline & Discrep. pr.  & $L$-curve &GCV& Quasi-opt
\\
\hline \hline
$\lambda$ & $ 6.3937\cdot 10^{-4} $& $ 4.9988\cdot 10^{-4} $&
$4.9988\cdot 10^{-4}$& $4.9988\cdot 10^{-4}$
\\
\hline $\frac{\|\Delta x_\lambda\|_2}{\epsilon \|x_\lambda\|_2}$
& $9.1464\cdot 10^{-1}$& $1.4820$& $1.1081$&  $3.1918$
\\
\hline $\frac{\|\Delta x_\lambda\|_\infty}{\epsilon \|x_\lambda\|_\infty}$
& $9.3703\cdot 10^{-1} $&  $1.6158$& $1.3755$& $3.7359$
\\
\hline  $\epsilon^{-1} \left\|\frac{\Delta x_\lambda}{x_\lambda}\right\|_\infty$
& $2.29 \cdot 10^6$ & $6.26 \cdot 10^6$ & $1.55 \cdot 10^6$ & $4.74 \cdot 10^6$
\\
\hline $\rm cond^F_{_{Reg}}$ in (\ref{condF})& $3.3961\cdot 10^3$& $4.4761\cdot 10^3$&
$4.4761\cdot 10^3$& $4.4761\cdot 10^3$
\\
\hline $\rm m_{_{Reg}}$ in
(\ref{mReg})& $1.5204\cdot 10^3$& $2.0035\cdot 10^3$& $2.0035\cdot 10^3$& $2.0035\cdot 10^3$
\\
\hline  $\rm c_{_{Reg}}$ in (\ref{cReg}) & $9.8192\cdot 10^6$&
$1.6064\cdot 10^7$& $1.6064\cdot 10^7$& $1.6064\cdot 10^7$
\\
 \hline
 $\kappa_{\sf SymToep}(A,b)$ & $1.0047\cdot 10^3$& $1.3242\cdot
10^3$&  $1.3242\cdot 10^3$&  $1.3242\cdot 10^3$
\\
\hline $m_{\sf SymToep}(A,b)$&  4.3765& 4.4971& 4.4971&4.4971
\\
\hline  $c_{\sf SymToep}(A,b)$  & $9.8143\cdot 10^6$&
$1.6056\cdot 10^7$& $1.6056\cdot 10^7$&$1.6056\cdot 10^7$
\\
 \hline
\end{tabular}
\end{table}

Table~\ref{Ex:toep} shows that the structured mixed condition numbers are
much smaller than the corresponding unstructured ones, which give
tight linear perturbation bounds. Both $c_{\sf SymToep}(A,b)$ and
$\kappa_{\sf SymToep}(A,b)$ are smaller than the corresponding
unstructured ones.

In Table~\ref{Ex:toep second}, we choose $M=e_3^\rt$, where $e_3$ is the
third column of $I_5$. In this case,
$$\frac{\|M\Delta
x_\lambda\|_\infty}{\|M x_\lambda\|_\infty}=\left\|\frac{M\Delta
x_\lambda}{M x_\lambda}\right \|_\infty, \qquad  m_{\sf
SymToep}(A,b)=c_{\sf SymToep}(A,b).$$
We compare the true relative
perturbation bounds with the first-order asymptotic perturbation
bounds given by $m_{\sf SymToep}(A,b)$ and $\kappa_{\sf SymToep}(A,b)$.

\begin{table}\centering
\caption{\label{Ex:toep second}  $L=I_5$, $M=e_3^\rt$.}
\begin{tabular}{||c|c c c c||}
\hline & Discrep. pr.  & $L$-curve &GCV& Quasi-opt
\\
\hline \hline
$\lambda$&  $6.39\cdot 10^{-4}$& $5.00\cdot 10^{-4}$ &
 $5.08\cdot 10^{-4}$&  $5.08\cdot 10^{-4}$
\\
\hline $\frac{\|M\Delta x_\lambda\|_2}{\|Mx_\lambda\|_2}$
& $1.3741\cdot 10^{-2}$& $2.2481\cdot 10^{-2}$& $2.1749\cdot
10^{-2}$& $2.1749\cdot 10^{-2}$
\\
\hline $\frac{\|M\Delta x_\lambda\|_\infty}{\|M x_\lambda\|_\infty}$
& $1.3741\cdot 10^{-2} $& $2.2481\cdot 10^{-2}$& $2.1749\cdot
10^{-2}$& $2.1749\cdot 10^{-2}$
\\
\hline
$\varepsilon\kappa_{\sf SymToep}(A,b)$ & $1.7780\cdot
10^{-1}$& $2.9088\cdot 10^{-1}$  & $2.8141\cdot 10^{-1} $&
$2.8141\cdot 10^{-1}$
\\
\hline $\varepsilon m_{\sf SymToep}(A,b)$& $ 4.9096\cdot 10^{-2 }$&
$8.0320\cdot 10^{-2}$& $ 7.7705\cdot 10^{-2}$& $7.7705\cdot 10^{-2}$
\\
 \hline
\end{tabular}
\end{table}

From Table \ref{Ex:toep second}, we can see the quantities
$\varepsilon m_{\sf SymToep}(A,b)$ give tighter perturbation bounds than
$\varepsilon\kappa_{\sf SymToep}(A,b)$, since they have the same
order as that of the true relative perturbation bounds.

Table~\ref{Ex:toep third} shows the results from different
choices of $M$, i.e., $M=e_1^\rt$
and $M=e_3^\rt$. For example, if we choose $M=e_1^\rt$, then we are
interested in the conditioning  of the first component of
$x_\lambda$. We display the values of $m_{\sf SymToep}(A,b)$ and
$\kappa_{\sf SymToep}(A,b)$. From Table 3, we can say that the first
component of $x_\lambda$ has better conditioning than the third one.

\begin{table}\centering
\caption{\label{Ex:toep third}  $L=I_5$, $M=e_1^\rt$ and $M=e_3^\rt$.}
\begin{tabular}{||c|c c c c||}
\hline $M=e_1^\rt $ & Discrep. pr.  & $L$-curve &GCV& Quasi-opt
\\
\hline \hline
$\lambda$& $  6.39\cdot 10^{-4 }$& $5.00\cdot 10^{-4}$ &
 $5.08\cdot 10^{-4} $&  $5.08\cdot 10^{-4}$
\\
\hline
$\kappa_{\sf SymToep}(A,b)$ & $1.5887\cdot 10^3$ &
$2.0941\cdot 10^3$ & $2.0594\cdot 10^3$ & $2.0594\cdot 10^3$
\\
\hline $m_{\sf SymToep}(A,b)$& $7.6056\cdot 10^2$&
$1.0022\cdot 10^3$ & $9.8567 \cdot 10^2$&  $9.8567\cdot 10^2$
\\
\hline
$M=e_3^\rt $ & Discrep. pr.  & $L$-curve &GCV& Quasi-opt
\\
\hline Value of $\lambda$&  $6.39\cdot 10^{-4}$& $5.00\cdot 10^{-4}$ &
$5.08\cdot 10^{-4}$&  $5.08\cdot 10^{-4}$
\\
\hline
$\kappa_{\sf SymToep}(A,b)$ & $1.7780\cdot 10^7$&
$2.9088\cdot 10^7 $& $2.8141\cdot 10^7$& $2.8141\cdot 10^7$
\\
\hline $m_{\sf SymToep}(A,b)$&  $4.9096\cdot 10^6$&
$8.0320\cdot 10^6 $&  $7.7705\cdot 10^6$& $7.7705\cdot 10^6$
\\
\hline
\end{tabular}
\end{table}

At the end of this example, we use Algorithms \ref{algorithm:1} and
\ref{al:mixed} to illustrate
the effectiveness of the power method. We set the maximal number of
iterations to 10 in Algorithm \ref{algorithm:1}. The estimated
condition numbers in Algorithms \ref{algorithm:1}
and \ref{al:mixed} are denoted by $\kappa_{\sf SymToep}^{\rm Est}(A,b)$ and
$m_{\sf SymToep}^{\rm Est}(A,b)$ respectively.

From Table \ref{Ex:toep estimation}, we can say that
$\kappa_{\sf SymToep}^{\rm Est}(A,b)$ and
$m_{\sf SymToep}^{\rm Est}(A,b)$ give good estimations
for this specific $A$, $L$ and $b$, especially
$m_{\sf SymToep}^{\rm Est}(A,b)$ gives better estimation.

\begin{table}\centering
\caption{\label{Ex:toep estimation}  $L=I_5$, $M=I_5$ and $M=e_1^\rt$.}
\begin{tabular}{||c|c c c c||}
\hline $M=I_5 $ & Discrep. pr.  & $L$-curve &GCV& Quasi-opt
\\
\hline \hline
$\lambda$&   $6.39\cdot 10^{-4}$& $5\cdot 10^{-4}$ &
$5.08\cdot 10^{-4} $&  $5.08\cdot 10^{-4}$
\\
 \hline
$\kappa_{\sf SymToep}(A,b)$ &$1.5878783796\cdot 10^3$&
$2.0931466345\cdot 10^3$  & $2.0584876249 \cdot 10^3$&
$2.0584876249\cdot 10^3$
\\
\hline
$\kappa_{\sf SymToep}^{\rm Est}(A,b)$ &$7.5891802517\cdot 10^2$&
$1.0002500616\cdot 10^3$ & $9.8369583257\cdot 10^2$ &
$9.8369583257\cdot 10^2$
\\
\hline $m_{\sf SymToep}(A,b)$&$7.6117517197\cdot 10^2$&
$1.0027483753\cdot 10^3$ & $9.8617760333\cdot 10^2$&
$9.8617760333\cdot 10^2$
\\
\hline $m_{\sf SymToep}^{\rm Est}(A,b)$& $7.6055529470\cdot 10^2$ &
$1.0022493745\cdot 10^3$&  $9.8567031098\cdot 10^2$&
$9.8567031098\cdot 10^2$
\\
\hline
$M=e_1^\rt $ & Discrep. pr.  & $L$-curve &GCV& Quasi-opt
\\
\hline
$\lambda$&  $6.39\cdot 10^{-4}$& $5\cdot 10^{-4}$ &
$5.08 \cdot 10^{-4}$ & $5.08 \cdot 10^{-4}$
\\
 \hline
 $\kappa_{\sf SymToep}(A,b)$ & $1.5886924101\cdot 10^3$&
$2.0940875560\cdot 10^3$  & $2.0594198505\cdot 10^3$&
$2.0594198505\cdot 10^3$ 
\\
\hline
$\kappa_{\sf SymToep}^{\rm Est}(A,b)$ & $3.7980426062\cdot 10^2$&
$5.0050048168\cdot 10^2$   & $4.9222129601\cdot 10^2$&
$4.9222129601\cdot 10^2$  
\\
\hline $m_{\sf SymToep}(A,b)$& $7.6055560483\cdot 10^2$&
$1.0022496243\cdot 10^3$&   $9.8567056493\cdot 10^2$ &
$9.8567056493\cdot 10^2$
\\
\hline $m_{\sf SymToep}^{\rm Est}(A,b)$& $7.6055529470\cdot 10^2$ &
$1.0022493745\cdot 10^3$  &    $9.8567031098\cdot 10^2$ &
$9.8567031098\cdot 10^2$
\\
\hline
\end{tabular}
\end{table}
\end{example}

\begin{example}~({\bf \cite{Rump03b}})
Let $A=g([c;r])$ be the $6 \times 6$ Hankel matrix defined by
$$
  A=g([c;r])=\begin{bmatrix}
    h&1&1&-1&0&0\cr
    1&1&-1&0&0&0\\
    1&-1&0&0&0&-1\\
    -1&0&0&0&-1&1\\
    0&0&0&-1&1&1\\
    0&0&-1&1&1&0
  \end{bmatrix} , \quad
  c=\left[\begin{matrix}h\\ 1\\1\\-1\\ 0\\0\end{matrix}\right], \
  r=\begin{bmatrix}0\\0\\-1\\1\\1\\0\end{bmatrix},\
  b=\begin{bmatrix}h\\2\\0\\0\\2\\0\end{bmatrix}
$$
for $h=10^{-3}$, where $c$ is the first column of $A$ and $r$ is the
last row of $A$.

We can choose the basis
$Y_1=g([e_1;{\bf 0}), \ldots, Y_5 = g([e_5 ; {\bf 0}])$,
$Y_{6} = g ([e_6;e_1])$, $Y_7 = g([{\bf 0}; e_2]), \ldots$, $Y_{11}
 = g([{\bf 0};e_6])$,
so that
$
  A=\sum_{k=1}^{11}a_kY_k.
$
 Again, from Theorem \ref{theorem LS: linea struct} with
 the $Y_i$, the expressions for
 $m_{\hankel}(A,b),\, c_{\hankel}(A,b)$ and $\kappa_{\hankel}(A,b)$ can be obtained.

From Table~\ref{Ex:hankel}, we conclude that the
structured mixed condition numbers
can be much smaller than the corresponding unstructured condition
numbers. Structured mixed condition numbers also give sharp
perturbation bounds. The forward errors obtained by multiplying
the structured mixed
condition numbers with $10^{-8}$ are of the same order as
that of the exact errors.

\begin{table}\centering
\caption{\label{Ex:hankel}  $L=I_6$, $M=I_6$.}
\begin{tabular}{||c|c c c c||}
\hline & Discrep. pr.  & $L$-curve &GCV& Quasi-opt
\\
\hline \hline
$\lambda$&  $7.5918 \cdot 10^{-4} $ & $2.5002 \cdot 10^{-4}$&
$2.5002 \cdot 10^{-4}$&  $0.0017$
\\
\hline $\frac{\|\Delta x_\lambda\|_2}{\epsilon \|x_\lambda\|_2}$ & $1.0902$
& $2.9510$ & $2.6014$&  $1.3163$
\\
\hline $\frac{\|\Delta x_\lambda\|_\infty}{\epsilon \|x_\lambda\|_\infty}$
& $1.3264$ & $4.2237$& $3.2460$& $2.0419$
\\
\hline
$\epsilon^{-1} \left\|\frac{\Delta x_\lambda}{x_\lambda}\right\|_\infty$ &
$4.39 \cdot 10^6$ & $2.516 \cdot 10^7$ & $2.645 \cdot 10^7$ &
$1.31 \cdot 10^6$
\\
\hline
$\rm cond^F_{_{Reg}}$ in (\ref{condF})& $2.2310 \cdot 10^{3}$&
$1.1401 \cdot 10^{4}$&  $1.1401 \cdot 10^{4}$ & $4.6222 \cdot 10^2$
\\
\hline
$\rm m_{_{Reg}}$ in (\ref{mReg})& $7.8426 \cdot 10^2$ &
$4.0032 \cdot 10^{3}$& $4.0032 \cdot 10^{3}$& $1.6347 \cdot 10^2$
\\
\hline  $\rm c_{_{Reg}}$ in (\ref{cReg}) &
$1.3230 \cdot 10^{7}$& $1.0238\cdot 10^{8} $& $1.0238 \cdot 10^{8}$&
$2.6208 \cdot 10^{6}$
\\
\hline
$\kappa_{\hankel}(A,b)$ & $1.0372 \cdot 10^{3}$&
$5.2922 \cdot 10^{3}$&  $5.2922 \cdot 10^{3}$& $2.1648 \cdot 10^2$
\\
\hline $m_{\hankel}(A,b)$& 3.4999& 5.1247&5.1247&3.5000
\\
\hline  $c_{\hankel}(A,b)$  & $1.1576 \cdot 10^{7}$&
$8.9578 \cdot 10^{7}$& $8.9578 \cdot 10^{7}$& $2.2931 \cdot 10^{6}$
\\
\hline
\end{tabular}
\end{table}

\end{example}

\begin{example}{\bf (\cite{Vander})}\label{ex1}
Let $V=g(a)$ be a $25\times 10$ Vandermonde matrix whose
$(i,j)$-entry is
\[
V_{ij}=\left(\frac{j}{10}\right)^{i-1}, \
a=\left[\frac{1}{10},\frac{2}{10},\ldots,\frac{9}{10}, 1\right]^\rt, \
b \in \R^{25} \mbox{ with } b_{2k-1}=-1, b_{2k}=1.
\]

In Table \ref{Ex:Vand}, when $\lambda$ is small, the problem is
ill-conditioned under unstructured perturbations. The structured
condition  numbers are much smaller than the unstructured ones. The
perturbation bounds given by the structured condition numbers coincide
with the relative errors from the columns Discrep. pr. and $L$-curve.
When we
use GCV and Quasi-opt to compute the regularization parameters
$\lambda$, which is equal to 5.69 in this example, the problem is
well-conditioned. The structured condition
numbers have the same order as the unstructured ones.
Both of them give sharp perturbation bounds.

\begin{table}\centering
\caption{\label{Ex:Vand}  $L=I_{10}$, $M=I_{10}$.}
\begin{tabular}{||c|c c c c||}
\hline & Discrep. pr.  & $L$-curve &GCV& Quasi-opt
\\
\hline \hline
$\lambda$& $1.36 \cdot 10^{-5}$
& $6.31 \cdot 10^{-5}$ & $5.69$&  $5.69$
\\
\hline
$\frac{\|\Delta x_\lambda\|_2}{\epsilon \|x_\lambda\|_2}$ &
$2.7484$ & $3.007$ & $1.4131$& $1.4131$
\\
\hline
$\frac{\|\Delta x_\lambda\|_\infty}{\epsilon \|x_\lambda\|_\infty}$
& $3.4057$ & $2.6762$ & $2.7502$ & $2.7502$
\\
\hline
$\epsilon^{-1} \left\|\frac{\Delta x_\lambda}{x_\lambda}\right\|_\infty$ &
$2.2796 \cdot 10$&  $2.1046 \cdot 10$& $5.3054$& $5.3054$
\\
\hline
$\rm cond^F_{_{Reg}}$ in (\ref{condF})& $4.8637 \cdot 10^{6}$&
$1.7445 \cdot 10^{6}$ & $2.6028 \cdot 10$& $2.6028 \cdot 10$
\\
\hline
$\rm m_{_{Reg}}$ in (\ref{mReg})& $5.5816 \cdot 10^{4}$&
$2.3085 \cdot 10^{4}$ & $3.9279 \cdot 10$&  $3.9279 \cdot 10$
\\
\hline
$\rm c_{_{Reg}}$ in (\ref{cReg}) & $5.0645 \cdot 10^{5} $&
$6.3061 \cdot 10^{4} $& $8.5328 \cdot 10 $& $8.5328 \cdot 10$
\\
 \hline
$\kappa_\VdM(A,b)$ & $4.7123 \cdot 10$& $5.2721 \cdot 10$&
$1.2499 \cdot 10$& $1.2499 \cdot 10$
\\
\hline
$m_{\VdM}(A,b)$& $1.4219 \cdot 10$& $1.4076 \cdot 10$&
$2.0828 \cdot 10$& $2.0828 \cdot 10$
\\
\hline
$c_{\VdM}(A,b)$  & $1.8428 \cdot 10^2$& $6.1557 \cdot 10$&
$4.0178 \cdot 10$& $4.0178 \cdot 10$
\\
\hline
\end{tabular}
\end{table}

\end{example}

\begin{example} {\bf (\cite{CD07})}
Let $A=g(a)$ be  a $10\times 8$ Cauchy matrix whose $(i,j)$-entry is
\begin{eqnarray*}
a_{ij}&=&\frac{1}{i+j-1},\quad
a=[u_1,\ldots,u_{10},v_{1},v_2,\ldots,v_8]^\rt,\,\mbox{ with }
  u_i=i,\,v_j=1-j,\\
  b&=&[1,-1,1,-1,1,-1,1,-1,1,-1]^\rt.
\end{eqnarray*}
Then $A$ is a rectangular Hilbert matrix.

From the second and third columns (Discrep. pr. and $L$-curve)
of Table \ref{Ex:Cauchy}, we can see that the structured
condition numbers are much smaller than the unstructured one and
they give sharp perturbation bounds. The first-order unstructured
asymptotic perturbation bounds severely overestimate the true
relative errors in both normwise and componentwise cases for the
numerical examples of the discrepancy principle and $L$-curve
methods. As for the last
two columns, since the regularization parameter $\lambda$ is large,
the problems are well-conditioned. The structured condition numbers
are of the same order as that of the unstructured ones.

\begin{table}\centering
\caption{\label{Ex:Cauchy}  $L=I_{8}$, $M=I_{8}$.}
\begin{tabular}{||c|c c c c||}
\hline & Discrep. pr.  & $L$-curve &GCV& Quasi-opt
\\
\hline \hline
$\lambda$& $2.46 \cdot 10^{-10}$ & $6.97\cdot 10^{-7}$  & $1.72$ &
$1.72$
\\
\hline
$\frac{\|\Delta x_\lambda\|_2}{\epsilon \|x_\lambda\|_2}$ &
$2.7472$& $5.5724$& $3.6995$ & $3.6995$
\\
\hline
$\frac{\|\Delta x_\lambda\|_\infty}{\epsilon \|x_\lambda\|_\infty}$ &
$2.5007$ & $6.3752$ & $2.2306$ & $2.2306$
\\
\hline
$\epsilon^{-1} \left\|\frac{\Delta x_\lambda}{x_\lambda}\right\|_\infty$ &
$1.0489 \cdot 10$ & $1.4879\cdot 10$ & $1.4934\cdot 10^{2}$&
$1.4934\cdot 10^{2}$
\\
\hline
$\rm cond^F_{_{Reg}}$ in (\ref{condF})& $2.9150 \cdot 10^{8} $ &
$4.5472 \cdot 10^{7} $& $2.7426 \cdot 10$ &   $2.7426\cdot 10$
\\
\hline
$\rm m_{_{Reg}}$ in (\ref{mReg})& $2.8775 \cdot 10^{7}$ &
$8.0534\cdot 10^{6}$ & $1.0465\cdot 10$& $1.0465 \cdot 10$
\\
\hline
$\rm c_{_{Reg}}$ in (\ref{cReg}) & $1.0584 \cdot 10^{8}$ &
$4.6471 \cdot 10^{7}$ &$8.7045\cdot 10^{2}$ & $8.7045\cdot 10^{2}$
\\
\hline
$\kappa_\Cauchy(A,b)$ & $3.8644 \cdot 10$ & $4.1630 \cdot 10$ &
$5.4144 \cdot 10$&  $5.4144 \cdot 10$
\\
\hline
$m_{\Cauchy}(A,b)$& $1.8131 \cdot 10$ & $2.2502 \cdot 10$&
$7.2573$ &   $7.2573$
\\
\hline
$c_{\Cauchy}(A,b)$  & $3.4802 \cdot 10^2$ & $8.0663 \cdot 10$&
$3.9086 \cdot 10^2$ &  $3.9086 \cdot 10^2$
\\
\hline
\end{tabular}
\end{table}
\end{example}

In the rest of this section, we will show our test results on
the proposed SCE algorithms for the conditioning estimation of
the Tikhonov regularization solution.  Both the unstructured
and structured Tikhonov regularization cases are considered.
For the unstructured Tikhonov regularization, the test problems originally come from discretization of
Fredholm integral equations of the first kind, and they lead to
discrete ill-posed problems. We use the test problems included in
the REGULARIZATION TOOLS package. In these numerical experiments,
a discrete ill-posed problem $Ax=b$ using one of
the many built-in test problems is firstly generated;
then the white noise is added to the
right-hand side with a perturbation $e$ whose elements are
normally distributed with zero mean and standard deviation chosen
such that the noise-to-signal ratio
$\|e\|_2 / \|b\|_2 =10^{-4}$, thus producing a more
`realistic' problem.

We generated the perturbations $\Delta A= \varepsilon
\times (E\odot A)$ and $\Delta b= \varepsilon \times (f \odot b)$,
where $\varepsilon=10^{-8}$, $E$ and $f$ are random matrices
whose entries are uniformly distributed in the open interval $(-1,1)$.

To measure the effectiveness of the estimators, we define
the over-estimation ratios
$$
r_\kappa:=\frac{\kappa_{\rm SCE}^{(k)}\cdot \varepsilon}
{\|\Delta x_\lambda\|_2/\|x_\lambda\|_2},\quad
r_m:=\frac{m_{\rm SCE}^{(k)}\cdot \varepsilon}
{\|\Delta x_\lambda\|_\infty/\| x_\lambda\|_\infty},\quad
r_c:=\frac{c_{\rm SCE}^{(k)} \cdot \varepsilon}
{\|\Delta x_\lambda/x_\lambda \|_\infty},
$$
where $k$ is the subspace dimension in Algorithms \ref{algo:subnorm} and
\ref{algo:subcomp}, $\kappa_{\rm SCE}^{(k)},\, m_{\rm SCE}^{(k)}$ and
$c_{\rm SCE}^{(k)}$ are the outputs from Algorithms \ref{algo:subnorm}
and \ref{algo:subcomp}. Typically the ratios
in $(0.1, ~10)$ are acceptable \cite[Chapter 19]{HighamBook}.

For unstructured Tikhonov regularization problems, we test the SCE for
several classical ill-posed problems included in the
REGULARIZATION TOOLS package: \textsf{deriv2},
\textsf{shaw} and \textsf{wing}.
Those three examples give square coefficient matrices $A$ and
right-hand side vectors $b$. For the matrix $L$ in \eqref{TikReg1},
we chose the identity matrix and
$$
L_1=\begin{bmatrix}1
&-1&&\cr&\ddots&\ddots&\cr&&1&-1\end{bmatrix} \in \R^{(n-1)\times n}
$$
which approximates the first derivative operator.
We adopted the following four values of the regularization parameter
$\lambda$:
$$
0.1, \quad 6\cdot10^{-2},\quad 1.7\cdot10^{-3},\quad 1.7\cdot10^{-4}.
$$
In Table \ref{Ex:SCE unstr}, we report the numerical results on
the ratios $ r_\kappa$, $ r_m$ and $ r_c$ for
examples with various dimensions and choices of $L$. The table shows
that the mixed  condition estimation $m_{\rm SCE}^{(k)}$ reflects
the true error bound accurately, while the componentwise condition
estimation $c_{\rm SCE}^{(k)}$ gives accurate error bounds for
most cases and the normwise condition estimation fails to reflect
the true error bound accurately. Specifically, $r_m$ are between
$0.94$ and $16.33$, implying that the condition estimation
$m_{\rm SCE}^{(k)}$ can be considered reliable \cite{HighamBook}.
The values of the componentwise condition estimation $c_{\rm SCE}^{(k)}$
are within $(16.54,~88.19)$  except for the case  \textsf{shaw},
where $n=512$, $L=I_n$, $k=3$ with all choices of $\lambda$,
the case \textsf{wing}, where $n=256$, $L=I_n$, $k=3$,
$\lambda=0.1,\, 6\cdot10^{-2}$ and the case \textsf{shaw}, where
$n=256$, $L=L_1$, $k=5$, $\lambda=1.7\cdot10^{-3}$, indicating
that the componentwise condition estimation $c_{\rm SCE}^{(k)}$
is effective for most cases. For the normwise condition estimation
$\kappa_{\rm SCE}^{(k)}$, most of the values of $r_\kappa$ are of
order $\Oh(10^2)$, and even some of them are of order $\Oh(10^3)$,
showing that The normwise condition estimation overly estimates
for most of the cases.

{\tiny
\begin{table}\centering
\caption{\label{Ex:SCE unstr}  SCE for the Tikhonov regularization problem.}
\begin{tabular}{||c|c |c |c ||}
\hline
 & $ r_\kappa$ & $ r_m$  & $ r_c$ \\
\hline
\textsf{deriv2}, $n=64,\, L=I_n,\, k=5$& &  &
\\
\hline
$\lambda=0.1$ & $7.7533 \cdot 10^2 $ & $2.0847 $ & $5.1984\cdot 10$
\\
\hline
$\lambda=6\cdot10^{-2}$ & $7.6234\cdot 10^2  $ & $ 2.1284 $&
$4.7450 \cdot 10$
\\
\hline
$\lambda=1.7\cdot10^{-3}$ & $ 6.7427\cdot 10^2 $ & $1.2970  $&
$2.4927 \cdot 10$
\\
\hline
$\lambda=1.7\cdot10^{-4} $& $9.2294 \cdot 10^{2}$& $1.1437$&
$3.9081\cdot 10$
\\
\hline
\textsf{deriv2}, $n=64,\, L=L_1,\, k=5$& &  &
\\
\hline
$\lambda=0.1$ & $  9.2479\cdot 10^2 $ & $2.2288 $ & $3.8756\cdot 10$
\\
\hline
$\lambda=6\cdot10^{-2}$ & $ 9.2057 \cdot 10^2  $ & $ 2.5900$&
$ 3.6441 \cdot 10$
\\
\hline
$\lambda=1.7\cdot10^{-3}$ & $ 6.2048 \cdot 10^2 $ & $ 1.0458 $&
$  2.9864\cdot 10$
\\
\hline
$\lambda=1.7\cdot10^{-4} $& $  9.2686\cdot 10^{2}$& $1.5548$&
$ 3.4564\cdot 10$
\\
\hline
\textsf{wing}, $n=128,\, L=L_1,\, k=5$& &  & \\
\hline
$\lambda=0.1$ & $ 2.3753\cdot 10^3 $ & $1.6337\cdot 10 $ &
$1.6544\cdot 10$
\\
\hline
$\lambda=6\cdot10^{-2}$ & $1.4154 \cdot 10^3  $ & $ 7.5189$&
$1.5910 \cdot 10$
\\
\hline
$\lambda=1.7\cdot10^{-3}$ & $1.3068 \cdot 10^3 $ & $1.6301 $&
$2.2538\cdot 10$
\\
\hline
$\lambda=1.7\cdot10^{-4}$& $9.6524 \cdot 10^{2}$& $1.3917 $&
$ 8.1316\cdot 10$
\\
\hline
\textsf{wing}, $n=256,\, L=I_n,\, k=3$& &  & \\
\hline
$\lambda=0.1$ & $  6.1606\cdot 10^2 $ & $1.0112 $ &
$2.6130\cdot 10^2$
\\
\hline
$\lambda=6\cdot10^{-2}$ & $5.4330\cdot 10^2  $ & $ 1.3470$&
$1.8994\cdot 10^2$
\\
\hline
$\lambda=1.7\cdot10^{-3}$ & $1.0126 \cdot 10^3 $ &
$ 1.8223$& $ 6.2581\cdot 10$
\\
\hline
$\lambda=1.7\cdot10^{-4} $& $1.0923 \cdot 10^{3}$& $2.1715$&
$7.2018\cdot 10$ \\
\hline
\textsf{shaw}, $n=512,\, L=I_n,\, k=3$& &  & \\
\hline
$\lambda=0.1$ & $  1.6632\cdot 10^2 $ & $2.4174$ &
$1.7467\cdot 10^{2}$
\\
\hline
$\lambda=6\cdot10^{-2}$ & $1.8554 \cdot 10^2  $ &
$ 2.8948 $& $ 3.7301\cdot 10^{2} $
\\
\hline
$\lambda=1.7\cdot10^{-3}$ & $1.0707 \cdot 10^2$ &
$ 3.1782$& $3.7127\cdot 10^2$
\\
\hline
$\lambda=1.7\cdot10^{-4}$& $2.9693\cdot 10^{2}$&
$1.6429$& $1.8405\cdot 10^2$ \\
\hline
\textsf{shaw}, $n=256,\, L=L_1,\, k=5$& &  & \\
\hline $\lambda=0.1$ & $ 2.8617\cdot 10$ & $2.9029$ &
$8.8186\cdot 10$
\\
\hline
$\lambda=6\cdot10^{-2}$ & $3.1563 \cdot 10$ &
$ 3.0425$& $8.3908 \cdot 10$
\\
\hline
$\lambda=1.7\cdot10^{-3}$ & $3.0292 \cdot 10$ & $1.0652$&
$2.0672\cdot 10^{2}$
\\
\hline
$\lambda=1.7\cdot10^{-4}$& $1.1406 \cdot 10^{2}$&
$9.4368\cdot 10^{-1}$& $3.0579\cdot 10$ \\
\hline
\end{tabular}
\end{table}
}

For structured Tikhonov regularization cases, we tested
the following Toeplitz matrix:
$$
A=(a_{i-j})\in \R^{m\times n},\, a_{i-j}=\rho^{|i-j|}.
$$
We used the right-hand side $b={\bf e} \in \R^{m}$ and $\rho=0.99999$.
This Toeplitz matrix is also a symmetric matrix. The Tikhonov regularization
parameter is determined by the four classical criteria. As discussed
in Subsection~\ref{sec:sub structure}, similar to
Algorithms \ref{algo:subnorm} and \ref{algo:subcomp},
we can use the SCE to obtain the structured normwise,
mixed and componentwise condition estimations denoted by
$\kappa_{{\sf SymToep},\rm SCE}^{(k)}$, $m_{{\sf SymToep},\rm SCE}^{(k)}$
and $c_{{\sf SymToep},\rm SCE}^{(k)}$ respectively.
The perturbations $\Delta a$ on $a$ and $\Delta b$ on $b$ were generated as
in \eqref{eq:sec numerical example perturbation}. As in the previous
example, let the overestimate ratios be defined by
$$
r_\kappa^{\sf SymToep}:=
\frac{\kappa_{{\sf SymToep},\rm SCE}^{(k)}\cdot \varepsilon}
{\|\Delta x_\lambda\|_2/\|x_\lambda\|_2},\quad
r_m^{\sf SymToep}:=
\frac{m_{{\sf SymToep},\rm SCE}^{(k)}\cdot \varepsilon}
{\|\Delta x_\lambda\|_\infty/\| x_\lambda\|_\infty},\quad
r_c^{\sf SymToep}:=
\frac{c_{{\sf SymToep},\rm SCE}^{(k)} \cdot \varepsilon}
{\|\Delta x_\lambda/x_\lambda \|_\infty},
$$
which measure the reliability of the condition estimators. In this example,
we always set $L=I_n$.

In Table \ref{Ex:SCE str}, except for two cases, the values of all
the three ratios are of order $\Oh(10)$, implying that the SCE structured
condition estimations are reliable.

{\tiny
\begin{table}\centering
\caption{\label{Ex:SCE str}
SCE for the structured Tikhonov regularization problem.}
\begin{tabular}{||c|c |c |c ||}
\hline
 & $ r_\kappa^{\sf SymToep}$ & $ r_m^{\sf SymToep}$  & $ r_c^{\sf SymToep}$ \\
\hline
$m=100,\,n=50,\, k=3$& &  &
\\
\hline
Discrep. pr. $\lambda=2.21$ & $1.3427 \cdot 10$ & $1.2354$ &
$4.0474\cdot 10$
\\
\hline
$L$-curve \, $\lambda=6.19 \cdot 10^{-2}$ & $1.8396 \cdot 10$
& $1.7149$& $2.9893 \cdot 10$
\\
\hline
GCV\, $\lambda=1.35\cdot 10^{-4}$ & $1.8113\cdot 10$ &
$1.5819$& $2.0690 \cdot 10$
\\
\hline
Quasi-opt\, $\lambda=7.48\cdot10^{-1}$& $1.1606 \cdot 10$&
$8.6429 \cdot 10^{-1}$& $5.7865\cdot 10$
\\
\hline
$m=300,\,n=200,\, k=3$& &  &
\\
\hline
Discrep. pr. $\lambda=4.71$ & $7.6537 \cdot 10$ & $4.1397$ &
$4.8346\cdot 10$
\\
\hline
$L$-curve  $\lambda=1.10 \cdot10^{-1}$ & $5.6455 \cdot 10$ &
$2.4765$& $7.8375 \cdot 10$
\\
\hline
GCV $\lambda=3.22\cdot 10^{-4 }$ & $3.2584\cdot 10$ & $1.4440$ &
$ 8.8256 \cdot 10$
\\
\hline
Quasi-opt\, $\lambda=4.49 $& $7.1689 \cdot 10$& $3.6997$&
$5.4007\cdot 10$
\\
\hline
$m=500,\,n=300,\, k=3$& &  &  \\
\hline
Discrep. pr. $\lambda=1.49\cdot 10^{-2}$ & $5.8732 \cdot 10$ &
$1.9115$ & $1.0526\cdot 10^2$
\\
\hline
$L$-curve $\lambda=1.03 $ & $7.8779 \cdot 10$ & $3.2648$&
$9.2462 \cdot 10$
\\
\hline
GCV $\lambda=5.66\cdot 10^{-4 }$ & $4.3191\cdot 10$ &
$1.6125$& $ 1.2357 \cdot 10$
\\
\hline
Quasi-opt $\lambda=9.25$& $1.0777 \cdot 10^2$& $4.3658$&
$6.8813\cdot 10$
\\
\hline
\end{tabular}
\end{table}

}

\section{Concluding Remarks}
In this paper, we introduce the structured condition numbers for the
structured Tikhonov regularization problem and derive their exact
expressions without the Kronecker product.
The structures considered include linear structures,
such as Toeplitz and Hankel,
and nonlinear structures, such as Vandermonde and Cauchy. We show that
our structured condition numbers are smaller than unstructured condition
numbers for
Toeplitz and Hankel structures. Applying the power method,
we devise fast algorithms for estimating the unstructured and structured
condition number under normwise and componentwise perturbations, that
can be integrated into a GSVD based Tikhonov regulariztion solver.
We also investigate the SCE for estimating structured condition
numbers. The numerical examples show that
our structured mixed condition numbers give tight error bounds and
the proposed condition estimations are reliable and efficient. A possible future research topic is to study the ratio between the structured and unstructured condition numbers for the structured Tikhonov regularization problem.

\end{document}